\documentclass[11pt,a4paper]{article}

\usepackage{geometry}
\usepackage{a4wide}

\usepackage[utf8]{inputenc}
\usepackage[numbers,square]{natbib}

\usepackage{amsmath,amsfonts,amsthm,amssymb} %For \eqref

\usepackage{hyperref}

\usepackage{booktabs}
\usepackage{tabularx}  
\usepackage{multirow}
\usepackage{longtable}
\usepackage{colortbl}

\usepackage[english]{babel}

\usepackage{calc}

\usepackage{algorithmic}
\usepackage{listings}				% zum Einbinden von Programmcode  
\RequirePackage{verbatim} 
\usepackage[%vlined,
algoruled,
]{algorithm2e}
\SetKwComment{tcp}{~\hspace{.8cm}~\textnormal{//} }{}%

\SetCommentSty{mycommfont}
\usepackage{rotating}
\LinesNumbered

\usepackage{tikz}		
\usepackage{graphicx}
\usepackage{subcaption}
\usepackage{wrapfig}
\usetikzlibrary{shapes.arrows,chains}
\usetikzlibrary{fit}
\usetikzlibrary{patterns}
\usetikzlibrary{shapes.geometric}
\usetikzlibrary{arrows,automata}
\usetikzlibrary{shapes,positioning}
\usetikzlibrary{plotmarks}

\usepackage{pgfplots}
\usepgfplotslibrary{fillbetween}
\pgfplotsset{compat=newest}

\newtheorem{myprop}{Proposition}[section]

\newtheorem{mydef}[myprop]{Definition}

\usepackage{accents}
\newcommand\munderbar[1]{\underaccent{\bar}{#1}}
\usepackage{authblk}

\title{Multistage Robust Discrete Optimization\\via Quantified Integer Programming}

\author[1]{Marc Goerigk}
\author[2]{Michael Hartisch\footnote{Corresponding author. Email: \texttt{michael.hartisch@uni-siegen.de}}}
\affil[1]{Network and Data Science Management, University of Siegen, Germany}
\affil[2]{Technology Management, University of Siegen, Germany}
\date{}

\usepackage[normalem]{ulem} %for sout

\newcommand{\QIPPU}{QIP$^{\text{PU}}$\xspace}
\newcommand{\QIPPUNew}{QIP$^{\text{PU}}$}

\newcommand{\PSPACE}{\mbox{\texttt{PSPACE}}\xspace}
\newcommand{\NP}{\mbox{\texttt{NP}}\xspace}

\newcommand{\EG}{\mbox{e.g.}\xspace}
\newcommand{\IE}{\mbox{i.e.}\xspace}

\newcommand{\X}{\mathcal{X}}
\newcommand{\cU}{\Xi}
\newcommand{\cR}{\mathcal{R}}

\def\D_All{\mathcal{D}}

\allowdisplaybreaks

\newcolumntype{R}[1]{>{\raggedleft\arraybackslash}p{#1}}
\newcolumntype{L}[1]{>{\raggedright\arraybackslash}p{#1}}

\makeatletter
\newcommand{\ssymbol}[1]{^{\@fnsymbol{#1}}}
\makeatother    
\newlength{\WidthStar}
\newlength{\WidthDagger}
\newcommand{\KNAQPU}{\textsc{\QIPPUNew(Kna)}\xspace}
\newcommand{\KNARC}{\textsc{DEP(Kna)}\xspace}

\newcommand{\LOTQ}{\textsc{QIP(Lot)}\xspace}
\newcommand{\LOTRC}{\textsc{DEP(Lot)}\xspace}

\newcommand{\SELQ}{\textsc{QIP(Sel)}\xspace}
\newcommand{\SELQPU}{\textsc{\QIPPUNew(Sel)}\xspace}
\newcommand{\SELRC}{\textsc{DEP(Sel)}\xspace}

\newcommand{\ASSQ}{\textsc{QIP(Ass)}\xspace}
\newcommand{\ASSQPU}{\textsc{\QIPPUNew(Ass)}\xspace}
\newcommand{\ASSRC}{\textsc{DEP(Ass)}\xspace}
\begin{document}

\maketitle
\begin{abstract}
Decision making needs to take an uncertain environment into account. Over the last decades, robust optimization has emerged as a preeminent method to produce solutions that are immunized against uncertainty. The main focus in robust discrete optimization has been on the analysis and solution of one- or two-stage problems, where the decision maker has limited options in reacting to additional knowledge gained after parts of the solution have been fixed. Due to its computational difficulty, multistage problems beyond two stages have received less attention.

In this paper we argue that multistage robust discrete problems can be seen through the lens of quantified integer programs, where powerful tools to reduce the search tree size have been developed. By formulating both integer and quantified integer programming formulations, it is possible to compare the performance of state-of-the-art solvers from both worlds. Using selection, assignment, lot-sizing and knapsack problems as a testbed, we show that problems with up to nine stages can be solved to optimality in reasonable time.
\end{abstract}


\section{Introduction}

Uncertainty affects most aspects of decision making, and thus needs to be taken into account preemptively. Different methodologies have been developed for this purpose, such as stochastic programming \cite{KM05} or robust optimization \cite{Ben-Tal}, which is the focus of this paper.
We consider discrete optimization problems of the form
\begin{align*}
\min_{\pmb{x}\in\X}\ &\pmb{c}(\pmb{\xi}) \pmb{x} \\
\text{s.t. } &A(\pmb{\xi}) \pmb{x} \le \pmb{b}(\pmb{\xi})
\end{align*}
%\[ \min_{\pmb{x}\in\X} \pmb{c} \pmb{x} \]
where $\X\subseteq\mathbb{Z}^n$ is the decision space, $\pmb{c}(\pmb{\xi})$ is a cost vector, $A(\pmb{\xi})$ the constraint matrix, and $\pmb{b}(\pmb{\xi})$ the right-hand side. In the following, vectors are always written in bold font and the transpose sign for the scalar product between vectors is dropped for ease of notation. The problem is affected by uncertainty, expressed by an uncertain parameter $\pmb{\xi}$.

The field of robust optimization incorporates a diverse set of approaches to formulate robust counterparts for such problems~\cite{goerigk2016algorithm}. In the most basic model we assume that a set of possible scenarios $\cU$ is given, the so-called uncertainty set. The (min-max, one-stage) robust counterpart is then to find a solution $\pmb{x}\in\X$ that performs best under its worst-case outcome, \IE to solve
\begin{align*}
\min_{\pmb{x}\in\X}\ &\max_{\pmb{\xi}\in\cU} \pmb{c}(\pmb{\xi}) \pmb{x} \\
\text{s.t. } &A(\pmb{\xi}) \pmb{x} \le \pmb{b}(\pmb{\xi}) & \forall \pmb{\xi}\in\cU
\end{align*}
%\[ \min_{\pmb{x}\in\X} \max_{\pmb{c}\in\cU} \pmb{c} \pmb{x} \]
This approach does not take a dynamic decision making process into account, where it is possible to react to new information when it becomes revealed. For this reason, two-stage approaches have been introduced \cite{yanikouglu2019survey}, which usually increase both the theoretical and practical complexity of solving robust problems. In many cases, it is possible to formulate such problems as mixed-integer programs (MIPs).

Different to the world of stochastic programming, where multistage problems beyond two stages are well-established~\cite{pflug2014multistage}, only few such approaches have been considered in the robust world. An intuitive reason is that the worst-case performance of a solution may depend on a single outcome of the scenario tree, and thus approximation methods such as sampling approaches
are not as powerful as in the stochastic world.

In this paper we show that finding optimal solutions to multistage robust discrete problems is actually well within the reach of current computational prowess. In particular, it is possible to see multistage robust problems through the lens of quantified integer programming (QIP),
which is an extension of integer programming where the variables are ordered explicitly and some variables are existentially and others are universally quantified.
QIPs are known to be $\PSPACE$-complete~\cite{Lorenz} and also can be interpreted as two-person zero-sum games between an existential and a universal (or adversarial) player.

The structure of this paper is as follows. In Section~\ref{sec:literature} we review related literature. We formally introduce QIPs in Section~\ref{sec:qip}, highlighting how to construct a QIP from a general multistage robust problem formulation. We illustrate this process using multistage selection, assignment, lot-sizing, and knapsack problems in  Section~\ref{sec:problems}. We introduce different QIP-based models and an extended MIP formulation for the robust counterpart. In Section~\ref{sec:exp} we discuss extensive computational experiments that compare CPLEX as a state-of-the-art general MIP solver with Yasol, a solver developed for QIPs. The paper is concluded in Section~\ref{sec:conclusion}.

\section{Related Literature}\label{sec:literature}

This section provides an overview of recent work in the areas mainly related to this paper: robust optimization (Section~\ref{sec:ro}) and quantified programming (Section~\ref{sec:quant}). For more general surveys on optimization under uncertainty, we refer to \cite{keith2019survey,ning2019optimization,bakker2019structuring}.

\subsection{Robust Optimization}\label{sec:ro}

Robust optimization problems are mathematical optimization problems with uncertain data, where a valid solution is sought for any (anticipated) realization of that data as represented by the uncertainty set $\cU$ \cite{ben2002robust}. 
Solving the robust counterpart ensures performance of the solution regarding $\cU$, but may result in a high price of robustness \cite{bertsimas2004price}, \IE the solution is often too conservative. Different concepts were developed to overcome this problem, \EG the concepts of light robustness \cite{fischetti2009light}, soft robustness \cite{ben2010soft}, adjustable robustness \cite{yanikouglu2019survey} and recoverable robustness \cite{Liebchen}.  In \cite{chassein2016performance} the authors discuss which approach is suitable for the problem at hand. A compact overview of prevailing uncertainty sets and robustness concepts can be found in \cite{gorissen2015practical,goerigk2016algorithm} where the latter focuses on the algorithm engineering methodology with regard to  robust optimization. Furthermore, in \cite{kasperski2016robust} a survey on robust discrete optimization is presented.

Two-stage models, \EG adjustable robust optimization and recoverable robust optimization, are often challenging to solve as even for simple cases the problem is $\NP$-hard \cite{ben2004adjustable}. Nevertheless, within the last few years several results regarding multistage models were obtained (\EG \cite{bertsimas2010finite,Ben-Tal}) and a tutorial-like survey on multistage robust decision-making was conducted in \cite{delage2015robust}.  A discussion of multistage optimization can be found in \cite[pp.~408--410]{Ben-Tal}
where the authors acknowledge its ``\textit{extreme applied importance}'' but point out the computational problems that arise and question the usefulness of most  approximation techniques. 
Besides frequently used variants of lot-sizing problems (\EG \cite{de2017robust,bertsimas2019adaptive})   multistage robust optimization has been applied to the daily operation of power systems \cite{lorca2016multistage}, resource allocation problems \cite{lin2008multi}, as well as planning and scheduling problems \cite{Lappas,ning2017data,motamed2019multistage}. An alternative approach in multistage discrete optimization is to construct a solution in the first stage, and modify it in further stages. Here one aims at finding stable solutions, which require little modification to remain near optimal for a changing cost function. Examples for this setting include matroids and matchings \cite{10.1007/978-3-662-43948-7_47}, the facility location problem \cite{eisenstat2014facility}, and the knapsack problem \cite{bampis_et_al:LIPIcs:2019:10966}. Dynamic programming techniques can be use to solve  multistage problems under uncertainty \cite{shapiro2011dynamic,mankowski2020extensions}, but often suffer from the curse of dimensionality. Other solution methods include variations of Bender's decomposition \cite{thiele2009robust}, column-and-constraint generation \cite{zeng2013solving}, and Fourier–Motzkin elimination \cite{zhen2018adjustable}. Additionally, iterative splitting of the uncertainty set is used to solve multistage robust problems in  \cite{postek2016multistage} and a partition-and-bound algorithm is presented in \cite{bertsimas2016multistage}. By considering specific robust counterparts a solution can be approximated and sometimes even guaranteed \cite{ben2004adjustable,
chen2009uncertain}.  Furthermore, several approximation schemes based on (affine) decision rules can be found in the literature (\EG \cite{bertsimas2015design,georghiou2019decision}).

\subsection{Quantified Programming}\label{sec:quant}
Quantified programs have been studied since at least 1995, when a first polynomial time algorithm for a restricted class of quantified linear programs based on quantifier elimination techniques  was introduced in \cite{Gerber}. The term \textit{quantified linear programming} (QLP) was coined in \cite{Subramani_Linear} and extended to quantified integer programming (QIP) in \cite{Subramani_Integer}. In  \cite{ESA11}  an objective function was introduced to the framework of quantified programming. Furthermore, a geometric analysis for QLPs \cite{Lorenz} 
and quantified programming models for games \cite{ederer2011modeling,lorenz2013solution} and combinatorial problems \cite{ederer2014quantified} have been presented. Algorithms for general QLPs were developed and implemented: an alpha–beta nested Bender's decomposition was proposed to solve the quantified linear optimization problem and tested in a computational study \cite{Euro}. A polyhedral uncertainty set was introduced for QIPs \cite{CG16}, 
which is closely related to the quantified linear implication problem \cite{eirinakis2012computational}.
Furthermore, an open-source solver for quantified programs was introduced in \cite{YasolACG17} and QIP specific pruning techniques were presented in \cite{hartisch2019novel}. 

Related to the concept of quantified programming is the quantified Boolean formula problem (QBF), which can be viewed as the satisfiability problem of a Boolean QIP where each constraint is a clause. Beside being the prototypical $\PSPACE$-complete problem \cite{stockmeyer1976polynomial} QBF allows a very compact problem description and thus several areas of application arise \cite{shukla2019survey}.

\section{Quantified Integer Programming}\label{sec:qip}

In the following, we formally introduce quantified integer programming. 
The thesis~\cite{DissMichael}, on which parts of this paper are based,
can be consulted for a more detailed discussion.

\subsection{Definition and Notation}

Let $n\in \mathbb{N}$ be the number of variables and $\pmb{x} \in \mathbb{Z}^n$ a vector of integer variables. We use the notation $[n]:=\{1,\ldots,n\}$ to denote index sets. 
For each variable $x_j$ its domain $\mathcal{L}_j$ with $l_j, u_j \in \mathbb{Z}$, $l_j \leq u_j$, $j\in[n]$, is given by
$\mathcal{L}_j =\{y\in \mathbb{Z} \mid    l_j \leq y \leq u_j \} \neq \emptyset$. 
The domain of the variable vector is  $\mathcal{L} =\{\pmb{y} \in \mathbb{Z}^n \mid \forall j \in [n]: y_j \in \mathcal{L}_j\}.$
 Let $\pmb{Q} \in \{\exists, \forall\}^n$ denote a vector of quantifiers.
We call each maximal consecutive subsequence in $\pmb{Q}$ consisting of identical quantifiers a \textit{quantifier block}. The quantifier corresponding to the $i$-th quantifier block is given by $Q^{(i)} \in \{\exists, \forall\}$  and the corresponding $i$-th \textit{variable block} is given by the (ordered) index set $B_i  \subseteq [n]$. 
Let $\beta \in [n]$ denote the number of variable blocks and thus $\beta-1$ is the number of quantifier changes. Note that $B_1\, {\cup}\, B_2\,  {\cup}\ldots  {\cup}\,B_\beta = [n]$ with $B_i \cap B_{i'} = \emptyset$ for $i \neq i'$.  
With $\mathcal{L}^{(i)}$ we denote the corresponding domain of the $i$-th variable block as in $\mathcal{L}$.

\begin{mydef}[Quantified Integer Linear Program (QIP)] \label{Def_QIP}~\\
Let $A^\exists \in \mathbb{Q}^{m_\exists \times n}$ and $\pmb{b}^\exists \in \mathbb{Q}^{m_\exists}$ for $m_\exists \in \mathbb{N}$. Let $\mathcal{L}$ and $\pmb{Q}$ be given with  $Q^{(1)}=Q^{(\beta)}=\exists$. Let $\pmb{c} \in \mathbb{Q}^n$ be the vector of objective coefficients, for which $\pmb{c}^{(i)}$ denotes the vector of coefficients belonging to variable block $B_i$. The term $\pmb{Q} \circ \pmb{x} \in \mathcal{L}$ with the component-wise binding operator $\circ$ denotes the \emph{quantification sequence} $Q^{(1)}\pmb{x}^{(1)} \in \mathcal{L}^{(1)}\ \ldots\ Q^{(\beta)} \pmb{x}^{(\beta)} \in\mathcal{L}^{(\beta)}$, such that every quantifier $Q^{(i)}$ binds the variables  $\pmb{x}^{(i)}$ of block $i$  ranging in their domain $\mathcal{L}^{(i)}$. We call $(A^\exists,\pmb{b}^\exists,\pmb{c},\mathcal{L},\pmb{Q})$ with
\[
	z = {\min\limits_{\pmb{x}^{(1)} \in \mathcal{L}^{(1)}}\left( \pmb{c}^{(1)}\pmb{x}^{(1)}+ \max\limits_{\pmb{x}^{(2)} \in \mathcal{L}^{(2)}} \left( \pmb{c}^{(2)}\pmb{x}^{(2)} + \min\limits_{\pmb{x}^{(3)} \in \mathcal{L}^{(3)}} \left( \pmb{c}^{(3)}\pmb{x}^{(3)} + \ldots \min\limits_{\pmb{x}^{(\beta)} \in \mathcal{L}^{(\beta)}} \pmb{c}^{(\beta)} \pmb{x}^{(\beta)}\right)\right)\right)} \] 
$$
\textnormal{s.t.}\ \pmb{Q} \circ \pmb{x} \in \mathcal{L}:\ A^\exists \pmb{x} \leq \pmb{b}^\exists
$$
	a \emph{quantified integer linear program} (QIP) with objective function.
\end{mydef}
We call $A^\exists \pmb{x} \leq \pmb{b}^\exists$ the \textit{(existential) constraint system} and $\mathcal{E}= \{i \in [\beta] \mid Q^{(i)}=\exists \}$ the set of existential variable blocks and $\mathcal{A}= \{i \in [\beta] \mid Q^{(i)}=\forall \}$ the set of universal variable blocks. Further, we call variable $x_j$ an existential (universal) variable if the corresponding quantifier $Q_j$ is $\exists$ ($\forall$).

A QIP can be interpreted as a two-person zero-sum game between an \emph{existential player} setting the existentially quantified variables and a \emph{universal player} setting the universally quantified variables with payoff $z$. The variables are set in consecutive order according to the variable sequence $x_1, \ldots, x_n$. We say that a player makes the move $\pmb{x}^{(i)} = \pmb{y}$, if she fixes the variable vector  $\pmb{x}^{(i)}$ of block $i$ to $\pmb{y} \in \mathcal{L}^{(i)}$. At each such move, the corresponding player knows the settings of $\pmb{x}^{(1)}, \ldots, \pmb{x}^{(i-1)}$ before taking her decision $\pmb{x}^{(i)}$. 

Each fixed vector $\pmb{x} \in \mathcal{L}$, that is, when the existential player has fixed the existential variables and the universal player has fixed the universal variables, is called \emph{a play}.  If $\pmb{x}$ satisfies the linear constraint system $A^\exists \pmb{x} \leq \pmb{b}^\exists$, the existential player pays $z=\pmb{c} \pmb{x}$ to the universal player. If $\pmb{x}$ does not satisfy $A^\exists \pmb{x} \leq \pmb{b}^\exists$, we say \emph{the existential player loses} and the payoff is $+\infty$.
Therefore, it is the existential player's primary goal to ensure the fulfillment of the constraint system, while the universal player tries to violate some constraints. If the existential player is able to ensure that all constraints are fulfilled he tries to minimize $\pmb{c}\pmb{x}$, whereas the universal player tries to maximize her payoff. A game tree can be used to represent the chronological order of all possible moves, given by the quantification sequence $\pmb{Q} \circ \pmb{x} \in \mathcal{L}$. The nodes in the game tree represent a partially assigned variable vector and the branches correspond to assignment of variables according to their variable domain. A solution is a so-called winning (existential) strategy, that defines how to react to each possible move by the universal player, in order to ensure  \(A^\exists \pmb{x} \le \pmb{b}^\exists \). Hence, a solution is a subtree of the game tree with an exponential number of leaves with respect to the number of universal variables and their domains. 
If there is more than one solution, the objective function aims for a certain (the ``best'') one, whereat the value of a strategy is defined via the worst-case payoff at its leaves.

\subsection{QIP with Polyhedral Uncertainty}
In optimization under uncertainty, uncertain events and variables must be modeled and selected cautiously as their influence might become too powerful. In the worst-case optimization framework provided by QIPs the modeler must ensure that the modeled worst case is not undesirable or too conservative, \IE the \textit{price of robustness} must be appropriate  \cite{bertsimas2004price}.
In  \cite{CG16} a second linear constraint system $A^\forall \pmb{x} \leq \pmb{b}^\forall$ was introduced that restricts the universal variables to a polytope resulting in the QIP with polyhedral uncertainty.

\begin{mydef}[QIP with Polyhedral Uncertainty (\QIPPU)]\label{Def::QIPPU}~\\
Let $m_\forall \in \mathbb{N}_0$, $\pmb{b}^\forall \in \mathbb{Q}^{m_\forall}$ and $A^\forall \in \mathbb{Q}^{m_\forall  \times n}$ with 
\begin{equation}\label{Eq::ExistZero}
A^\forall_{k,j}=0 \quad \forall \, k\in [m_\forall], \,  j \in [n] : Q_j=\exists\, .
\end{equation} 
Let $\D_All =\{\pmb{x}\in \mathcal{L} \mid A^\forall \pmb{x} \leq \pmb{b}^\forall\} \neq \emptyset$. The term $\pmb{Q} \circ \pmb{x} \in \D_All$ with the component-wise binding operator $\circ$ denotes the \emph{quantification sequence} $Q^{(1)}\pmb{x}^{(1)} \in {\D_All}^{(1)}\ Q^{(2)}\pmb{x}^{(2)}\in {\D_All}^{(2)}(\pmb{x}^{(1)})\ \ldots\ Q^{(\beta)} \pmb{x}^{(\beta)} \in {\D_All}^{(\beta)}(\pmb{x}^{(1)},\ldots,\pmb{x}^{(\beta-1)})$ such that every quantifier $Q^{(i)}$ binds the variables  $\pmb{x}^{(i)}$ of block $i$ ranging in their domain ${\D_All}^{(i)}(\pmb{x}^{(1)},\ldots,\pmb{x}^{(i-1)})$, with
$${\D_All}^{(i)}(\tilde{\pmb{x}}^{(1)},\ldots,\tilde{\pmb{x}}^{(i-1)})= \begin{cases}
\mathcal{L}^{(i)} &\text{if } i \in \mathcal{E}\\
\{\pmb{y} \in \mathcal{L}^{(i)} \mid \exists \pmb{x}=(\tilde{\pmb{x}}^{(1)},\ldots,\tilde{\pmb{x}}^{(i-1)},\pmb{y},{\pmb{x}}^{(i+1)},\ldots,{\pmb{x}}^{(\beta)}) \in \D_All \} 
&\text{if } i \in \mathcal{A}  \, .
\end{cases}$$ 
We call
\[z\ =\ \min\limits_{\pmb{x}^{(1)} \in  {\D_All}^{(1)}} \left(\pmb{c}^{(1)}\pmb{x}^{(1)}+ \max\limits_{\pmb{x}^{(2)} \in  {\D_All}^{(2)}} \left( \pmb{c}^{(2)}\pmb{x}^{(2)}  + \min\limits_{\pmb{x}^{(3)} \in  {\D_All}^{(3)}} \left(  \pmb{c}^{(3)}\pmb{x}^{(3)} + \ldots \min\limits_{\pmb{x}^{(\beta)} \in  {\D_All}^{(\beta)}}  \pmb{c}^{(\beta)}\pmb{x}^{(\beta)}\right) \right)\right)\]
$$
 \text{s.t.}  \ \pmb{Q} \circ \pmb{x} \in \D_All  : A^\exists \pmb{x}\leq \pmb{b}^\exists
$$
a \emph{QIP with polyhedral uncertainty} (\QIPPU) given by the tuple $(A^\exists, A^\forall, \pmb{b}^\exists, \pmb{b}^\forall, \pmb{c}, \mathcal{L}, \pmb{Q})$.
\end{mydef}

With Condition~\eqref{Eq::ExistZero} each entry of $A^\forall$ belonging to an existential variable is zero. Therefore, the \emph{universal constraint system} $A^\forall \pmb{x} \leq \pmb{b}^\forall$ restricts universal variables in such way that their range only depends on previous universal variables:
when assigning a value to the universal variable $x_i$, there must exist a series of future assignments for $x_{i+1}, \ldots, x_n$ such that the resulting vector $\pmb{x}$ fulfills  $A^\forall \pmb{x} \leq \pmb{b}^\forall$.  This means that a universal variable assignment must not make it impossible to satisfy the system $A^\forall \pmb{x} \leq \pmb{b}^\forall$. Hence, a correlation among the scenarios of different stages can be modelled. With $\D_All \neq \emptyset$ at least one fixation of universal variables fulfills the universal constraint system and therefore a universal strategy to fulfill $A^\forall \pmb{x} \leq \pmb{b}^\forall$ exists. This is similar to demanding a non-empty polyhedral uncertainty set, which is a common condition in robust optimization.

\subsection{The Open-Source QIP Solver Yasol}\label{SubSec::Yasol}
The open-source solver Yasol\footnote{Sources and further information regarding the solver can be found on \url{http://www.q-mip.org} (accessed 28 January, 2021).} is a search-based solver for QIPs \cite{YasolACG17}.  There are no other general QIP solvers that we know of. The heart of the search algorithm is an arithmetic linear constraint database together with an alpha-beta algorithm, which has been successfully used in gaming programs, \EG chess programs for many years \cite{KNUTH1975293,Hydra}. The aim of the search is to find a path from the root to the best valued leaf that can be reached, under the assumption of best play by both sides. The variable assignment corresponding to such a leaf is called \textit{principal variation}. The two values $\alpha$ and $\beta$ maintained during the search indicate the best known bounds for the universal and existential player, respectively, allowing to prune subtrees that cannot influence the optimal solution which would have to be traversed during standard minimax-search. The alpha-beta algorithm is particular advantageous if the principal variation is visited first, as the emerging values for $\alpha$ and $\beta$ more frequently lead to cutoffs.  Hence, the order in which the nodes are evaluated has a strong impact on its performance. 

In order to realize fast backjumps when encountering a contradiction---as typically performed in SAT- and QBF-solvers (\EG \cite{Qube,chen2001conflict})---the alpha-beta algorithm was extended as outlined in \cite{YasolACG17}: when a contradiction is detected a reason in form of a constraint is added to the constraint database and the search returns to the node where the found contradiction is no longer imminent. Yasol deals with constraint learning (conflict learning) on the so-called primal side as known from SAT- and QBF-solving (\EG \cite{marques2009conflict,giunchiglia2002learning}), as well as with constraint learning (cutting planes) on the dual side known from MIP (\EG \cite{ceria1998cutting}). Several other techniques from various research fields are implemented, \EG the killer heuristic \cite{Killer}, restart strategies \cite{biere2008adaptive}, strong branching \cite{achterberg2005branching}, and a QIP-specific pruning technique \cite{hartisch2019novel}.

Yasol is currently able to solve multistage quantified mixed integer programs with the following properties: a) The basic structure must be a quantified linear program, \IE linear (existential and universal) constraints and objective function, existentially or universally quantified variables, all variables are bounded from below and above and $Q_1=Q_n=\exists$.
b) Integer variables are allowed in all existential and universal variable blocks.
c) Continuous variables are allowed only in the closing existential variable block.
Due to the exponential size of a solution (strategy) the output of a feasible instance is the optimal assignment of the first (existential) variable block, as well as the optimal worst-case outcome, \IE the value of the optimal strategy.

Yasol makes intensive use of a linear programming solver in order to asses the quality of a branching variable (\EG \cite{achterberg2005branching}), check the satisfiability of the existential constraint system in the current subtree or for the generation of bounds. These tools are black-box used, while not exploiting the possible integer solving abilities of the foreign solver for standard QIPs. 
A universal constraint system is maintained in a second constraint database. The compliance of a universal variable block assignment with its corresponding domain $\D_All^{(i)}$---in general---requires the solution of an integer program. This, however, is not the case if there always exists a universal variable block assignment and if each assignment not obeying $\D_All^{(i)}$ immediately violates one of the universal constraints. This was motivated in \cite{hartisch2019OR} and more details can be found in \cite{DissMichael}. The instances of interest in this paper exhibit this characteristic and therefore no integer solving abilities of the foreign solver is exploited in any of the presented experiments.

\subsection{Linking Quantified Programming and Multistage Robust Optimization}\label{SubSec::ConnectionQIPMultistage}

To put the potential of quantified programming into the perspective of multistage robust discrete programming, we consider the general form as used in \cite{bertsimas2016multistage}:
\begin{subequations}\label{BertsimasModell}
\begin{align}
\min_{\pmb{x} \in \mathcal{X}} \max_{\pmb{\xi} \in \Xi}& \sum_{t\in[T]} \pmb{c}^t(\pmb{\xi})\cdot \pmb{x}^t(\pmb{\xi}^1,\ldots,\pmb{\xi}^{t-1}) \\
\textnormal{s.t.} & \sum_{t\in[T]} A^t(\pmb{\xi})\cdot \pmb{x}^t(\pmb{\xi}^1,\ldots,\pmb{\xi}^{t-1}) \leq \pmb{b}(\pmb{\xi})\qquad \forall \pmb{\xi}=(\pmb{\xi}^1,\ldots,\pmb{\xi}^T)\in \Xi
\end{align}
\end{subequations}

Here, $A^t(\pmb{\xi})$, $\pmb{b}(\pmb{\xi})$, and $\pmb{c}^t(\pmb{\xi})$ are assumed to be affine functions of $\pmb{\xi}$. We assume $\Xi$ to contain the integer lattice points within a bounded polyhedron given by a linear constraint system. 
For $t \in [T]$ we define
$$\Xi^t(\pmb{\hat{\xi}}^1,\ldots,\pmb{\hat{\xi}}^{t-1}) =\{\pmb{\hat{\xi}}^t \mid \exists \pmb{\xi}=(\pmb{\hat{\xi}}^1,\ldots,\pmb{\hat{\xi}}^{t-1},\pmb{\hat{\xi}}^t,\pmb{\xi}^{t+1},\ldots,\pmb{\xi}^{T})\in\Xi\}\, ,$$
\IE $\Xi^t$ is the domain of the uncertain parameters in stage $t$ and (potentially) depends on the realization $\pmb{\hat{\xi}}^1,\ldots,\pmb{\hat{\xi}}^{t-1}$ of uncertain parameters in previous stages.
We assume the domain $\mathcal{X}$ to be bounded and given by no additional constraints other than the variable bounds.  Let $\mathcal{L}^{(t)}$ be the integer lattice points in  the hypercube for which $\pmb{x}^t \in \mathcal{L}^{(t)}$ holds for all $\pmb{x} \in \mathcal{X}$, \IE $\mathcal{L}^{(t)}$ is the bounded domain of the variables $\pmb{x}^t$ of stage $t$. A  quantified program corresponding to the multistage robust problem is given as follows:
\begin{subequations}
\begin{align}
\min\ &  \sum_{t\in[T]} \pmb{c}^t(\pmb{\xi}) \pmb{x}^t \\
\textnormal{s.t.}\ &  
\exists\ \pmb{x}^1 \in \mathcal{L}^{(1)}\quad \forall\, \pmb{\xi}^1 \in \Xi^1\quad   \exists\ \pmb{x}^2 \in \mathcal{L}^{(2)} \ldots  \forall\, \pmb{\xi}^{T-1} \in \Xi^{T-1} \quad   \exists\ \pmb{x}^T \in \mathcal{L}^{(T)}  \quad \forall\, \pmb{\xi}^{T} \in \Xi^{T} :\nonumber\\
&\sum_{t\in[T]} A^t(\pmb{\xi}) \pmb{x}^t  \leq \pmb{b}(\pmb{\xi})
\end{align}
\end{subequations}

In this quantified program the decision variables $\pmb{x}^t$ at stage $t$ no longer explicitly depend on the realization of uncertain parameters $\pmb{\hat{\xi}}^1,\ldots,\pmb{\hat{\xi}}^{t-1}$. In particular, it is no longer necessary to have an entire set of variables representing the realization of $\pmb{x}^t$ for each realization of $\pmb{\hat{\xi}}^1,\ldots,\pmb{\hat{\xi}}^{t-1}$. Instead the quantification sequence demands that there must exist an assignment of $\pmb{x}^t$ for all scenarios. Hence, when using a quantified program, nonanticipativity is readily guaranteed. In order to obtain a QIP or \QIPPU the non-linearities in the objective and left-hand side must be removed. Due to the integrality of the existential variables $\pmb{x}^t$ and the universal variables $\pmb{\xi}$, we can linearize $\pmb{c}^t(\pmb{\xi}) \pmb{x}^t$ and $A^t(\pmb{\xi})\pmb{x}^t$ using standard linear programming techniques. Note that we omitted the min-max alternation in the objective and only stated the optimization sense for the (existential) decision variables. 

This transformation is not necessarily presented in order to be applicable, but to highlight the connection between the two different formulations of multistage problems. As we show in the following section, a quantified model can often be stated in a straight forward manner.  In fact, one of the strengths of a quantified model lies in the simplicity to state them: there is no need to consider nonanticipativity, significantly less variables and constraints are required as not every scenario subcourse has to be modeled, and rules for the anticipated uncertain variables can be modeled via a universal constraints system.
 In the following section we present different multistage optimization problems modeled as QIP and \QIPPU that exhibit uncertainty in the objective, left-hand side matrix or right-hand side vector.

\section{QIP Models for Multistage Robust Problems}\label{sec:problems}

In this section, we illustrate the general process of deriving both a QIP and a multistage robust problem formulation as outlined in Section~\ref{SubSec::ConnectionQIPMultistage} using four specific optimization problems. Most details are presented for the first such application, which is the multistage selection problem.

 \subsection{Multistage Robust Selection Problems}\label{SubSec::Selection}

We first examine the robust selection problem
%\[ \min_{\pmb{x}\in\X}\  \sum_{i\in[n]} c_i x_i  \]
%with $\mathcal{X} = \{ \pmb{x}\in\{0,1\}^n : \sum_{i\in[n]} x_i = p\}$, 
where $p$ out of $n$ items must be selected, such that the costs are minimized. 
A thorough overview of the one-stage robust selection problem can be found in \cite{CHASSEIN2018423}. In our multistage setting we assume that uncertainty is only present in the objective. In an initial (existential) decision stage a set of items can be selected for fixed costs $\pmb{c}^0$. Then, in a universal decision stage, one of $N \in \mathbb{N}$ cost scenario is revealed. %via the indicator variables $\pmb{q}^t$
In the subsequent existential decision stage further items can be selected and  variables $x_{i}^t$ indicate the selection of item $i$ in period $t$. Those two stages can be repeated iteratively several times. Let $T$ be the number of such iterations, \IE the number of universal decision stages. We call $T$ the  time horizon and we refer to $t\in [T]$ as a period. Hence, there will be $2T+1$ variable blocks: the initial existential variable block and the $T$ periods consisting of each a universal and existential variable block. 
The universal domain for each period $t \in [T]$ is given by
\[ \D_All_t = \left\{\pmb{q}^t \in \{0,1\}^N \mid \sum_{k\in[N]} q^t_k=1 \right\}\]
and $\pmb{q}^t \in \D_All_t$ is the vector indicating the selected scenario. %As before, let $N$ be the number of scenarios per period, \IE at each iteration, one of $N$ scenarios is revealed. 
The cost of item $i$ in scenario $k$ of period $t$ is given by $c_{i,k}^t$. 
This multistage selection problem under uncertainty can be modeled as a quantified program with a polyhedral uncertainty set as follows:
\begin{subequations}
\begin{align}
\text{\hspace*{-1cm}(\SELQPU)\hspace*{.5cm}} \min \  & \sum_{i\in[n]} c_{i}^0 x^0_i +  \sum_{t\in[T]} z_t \label{Eq::Ex_Mod_Sel1}\\
\textnormal{s.t.}\ &\rlap{$\exists \pmb{x}^0 \in \{0,1\}^n \quad \forall \pmb{q}^1 \in \D_All_1  \quad \exists \pmb{x}^1 \in \{0,1\}^n \quad \forall \pmb{q}^2 \in \D_All_2 \ \cdots $} \nonumber\\ 
& \rlap{$\cdots  \ \forall \pmb{q}^T \in \D_All_T  \quad \exists \pmb{x}^T \in \{0,1\}^n \ \exists \pmb{z} \in \mathbb{R}_+^T$:}~ \nonumber\\
&\sum_{i\in[n]} \sum_{t=0}^Tx_{i}^t = p  \label{Eq::Ex_Mod_Sel3}\\
& \sum_{t =0}^T x_{i}^t \leq 1 && \forall \, i\in[n]  \label{Eq::Ex_Mod_Sel4}\\
& \sum_{i\in[n]} c_{i,k}^t x_{i}^t  \leq z_t + M_k^t(1-q_{k}^t)\ \ &&  \forall  \, k \in [N], \, t \in [T] \label{Eq::Ex_Mod_Sel5}
\end{align}
\end{subequations}
The Objective~\eqref{Eq::Ex_Mod_Sel1} consists of the expenses from the initial decision stage with invariable costs and the expenses of subsequent periods in which the cost for each item depends on the selected scenario. Note that we omit the min/max alternation in the objective and only specify the optimization orientation for the existential variables. The first Constraint~\eqref{Eq::Ex_Mod_Sel3} demands that overall exactly $p$ items must be selected. Constraint~\eqref{Eq::Ex_Mod_Sel4} prevents that an item is selected more than once. Constraint~\eqref{Eq::Ex_Mod_Sel5} enforces the link between the selected scenario, selected items and resulting costs in each period $t$. If $M_k^t$ is selected appropriately for each potential scenario $k$ (\EG $M_k^t\geq \sum_{i\in[n]} c^t_{i,k}$), all but one of the Constraints \eqref{Eq::Ex_Mod_Sel5} are trivially fulfilled for a realization of $\pmb{q}^t \in \D_All_t$: if scenario $k$ is selected by the universal player ($q^t_k=1$) the corresponding costs $\pmb{c}_{\cdot,k}$ are decisive for the cost calculation. Those costs are bundled in the 
 auxiliary $z_t$ variables,  which we put the at the end of the quantification sequence, as Yasol can only deal with (existential) continuous variables in the last variable block.  Here, however, the cost variables $z_t$ also could be placed immediately after the corresponding selection in period $t$.  Additionally, when explicitly stating the model, one has to specify an upper bound on $z_t$, which can be easily computed by taking the cost vectors of the corresponding scenarios into account.

In order to build an equivalent QIP, \IE a model without constraints on the universal variables, we use an integer variable $\ell_t \in [N]$ in order to select one of $N$ scenarios in each period. This integer can then be transformed into existential indicator variables, resulting in the following  problem:
\begin{subequations}
\begin{align}
\text{\hspace*{-1cm}(\SELQ)\hspace*{.5cm}} \min \  & \sum_{i\in[n]} c_{i}^0 x^0_i +  \sum_{t\in[T]} z_t \label{Eq::Ex_Mod_QSel1}\\
\textnormal{s.t.}\ &\rlap{$\exists \pmb{x}^0 \in \{0,1\}^n \quad \forall \ell_1 \in [N] \quad \exists \pmb{q}^1 \in \{0,1\}^N \quad \exists \pmb{x}^1 \in \{0,1\}^n \ \cdots  $}~ \nonumber\\
&\rlap{$\cdots \  \forall \ell_T \in [N] \quad \exists \pmb{q}^T \in \{0,1\}^N \quad \exists \pmb{x}^T \in \{0,1\}^n \ \exists \pmb{z} \in \mathbb{R}_+^T$:}~  \nonumber\\
&\sum_{i\in[n]} \sum_{t=0}^T x_{i}^t = p  \label{Eq::Ex_Mod_QSel3}\\
& \sum_{t =0}^T x_{i}^t \leq 1 && \forall \, i\in[n]  \label{Eq::Ex_Mod_QSel4}\\
& \sum_{i \in[n]} c_{i,k}^t x_{i}^t  \leq z_t + M_k^t(1-q_{k}^t) \ \ &&  \forall \,  k \in [N],\, t \in [T]  \label{Eq::Ex_Mod_QSel5}\\
& \sum_{k\in[N]} q^t_k=1 && \forall  \, t\in [T]  \label{Eq::Ex_Mod_QSel6}\\
& \sum_{k\in[N]} k \cdot q^t_k=\ell_t && \forall \, t\in[T]  \label{Eq::Ex_Mod_QSel7}
\end{align}
\end{subequations}
The variables $\pmb{q}^t$, which were universal variables in \SELQPU, are now used as existential variables indicating the selected scenario. Constraints \eqref{Eq::Ex_Mod_QSel6} and \eqref{Eq::Ex_Mod_QSel7} ensure that the 
scenario number $\ell_t$ selected by the universal player is transformed into a corresponding assignment of $\pmb{q}^t$. Thus, the number of variables and constraints in \SELQ increased compared to \SELQPU.

We also provide a deterministic equivalent program (DEP), \IE an equivalent MIP, for which each possible scenario sequence must be listed explicitly (see model \eqref{BertsimasModell}).  The set containing all possible sequences of scenarios is $\cR= [N]^T$. For one such sequence $\pmb{r} \in \cR$ the scenario in period $t$ is $r_t$.  The entire sub-sequence  up to period $t$ is denoted by $\pmb{r}(t) \in   [N]^t$. For each item $i \in [n]$, period $t \in [T]$ and sequence $\pmb{r} \in \cR$, the  variable $x_i^{\pmb{r}(t)}$ indicates the decision of selecting item $i$ in period $t$ after the sub-sequence $\pmb{r}(t)$ of $\pmb{r}$ occurred. Note that $x_i^{r(t)}$ corresponds to writing $x_i(\pmb{\xi}_1,\ldots,\pmb{\xi}_t)$ in model \eqref{BertsimasModell}. This ensures the nonanticipativity property: even for different scenario sequences the selection decisions must be the same, as long as the sub-sequences are identical.

As an example, for $N=4$ and $T=6$ a possible sequence of scenarios is $\pmb{r}=(1,4,2,3,1,1)$. The sub-sequence until period $t=4$ is $\pmb{r}(4)=(1,4,2,3)$. The variable indicating whether item $i$ is selected after $4$ periods and the occurrence of this particular sub-sequence is denoted $x_i^{\pmb{r}(4)}=x_i^{(1,4,2,3)}$ and the scenario in period $4$ for this sequence is $r_4=3$. The cost of item $i$ in period $t=4$ does not depend on the entire sequence, but only on the occurred scenario and is given by $c_{i,r_t}^{t}=c_{i,3}^4$. For the sequence of scenarios $\hat{\pmb{r}}=(1,4,2,3,2,4)$ it holds  $\pmb{r}(4) = \hat{\pmb{r}}(4)$ and therefore, the variables $x_i^{\pmb{r}(4)}$ and $x_i^{\hat{\pmb{r}}(4)}$ are the same. 

The DEP of \SELQPU represents the robust counterpart to the problem and is given as follows.
\begin{subequations}
\begin{align}
\text{\hspace*{-1cm}(\SELRC)\hspace*{.5cm}}\min\ &\sum_{i\in[n]} c_i^0 x_i^0 + z \label{Eq::Ex_Mod_MIPSel1}\\
\textnormal{s.t. } & z \ge \sum_{i\in[n]} \sum_{t\in[T]} c^{t}_{i,r_t} x^{\pmb{r}(t)}_i && \forall \, \pmb{r} \in \cR \label{Eq::Ex_Mod_MIPSel2} \\
& \sum_{i\in[n]} \left( x_i^0 + \sum_{t\in[T]} x^{\pmb{r}(t)}_i \right)  = p && \forall \, \pmb{r} \in \cR \label{Eq::Ex_Mod_MIPSel3}\\
& x_i^0 + \sum_{t\in[T]} x^{\pmb{r}(t)}_i  \le 1 && \forall \, i\in[n],\, \pmb{r} \in \cR \label{Eq::Ex_Mod_MIPSel4}\\
& x_i^0 \in\{0,1\} && \forall \, i\in [n] \label{Eq::Ex_Mod_MIPSel5}\\
& x^{\pmb{r}(t)}_i \in\{0,1\}  && \forall  \,i\in [n],  \,\pmb{r} \in \cR, \, t \in [T] \label{Eq::Ex_Mod_MIPSel6}\\
&z \in\mathbb{R}
\end{align}
\end{subequations}
Constraint~\eqref{Eq::Ex_Mod_MIPSel2} ensures that the expenses from the worst-case scenario sequence appear in the objective function. Constraint~\eqref{Eq::Ex_Mod_MIPSel3} ensures for each scenario sequence that exactly $p$ items are selected in the end, whereas Constraint~\eqref{Eq::Ex_Mod_MIPSel4} ensures  that each item is selected at most once. Note that this formulation grows exponentially in $T$, whereas the quantified programs allow for a compact problem description.

 \subsection{Multistage Robust Assignment Problems}\label{SubSec::Assignment}

Another frequently considered combinatorial problem is the assignment problem: Given a complete bipartite graph $G = (V, E)$ with  $V = A \cup B$, $n=|A| = |B|$. A cost value $c_{i,j} \in \mathbb{R}_+$ is associated with each edge $(i,j) \in E$. The assignment problem consists of determining a perfect matching of minimum costs. Research on the min-max and min-max regret assignment problems can be found in
\cite{aissi2005complexity} and further complexity results are obtained in \cite{dei2006robust}, while recently, also recoverable variants have been considered~\cite{bold2020recoverable,fischer2020investigation}. We want to adapt the robust approach to a multistage setting. 

In an initial decision stage the existential player can select edges with known costs. Then, iteratively, new costs of the edges are presented (by the universal player) which then can be selected (by the existential player). Similar to the preceding subsection, the costs selected by the universal player come from a predefined scenario pool. Let $N$ be the number of scenarios and $T$ the time horizon. Universal variables $q^t_k$ indicate whether cost scenario $k$ is selected in period $t$. As only one scenario can occur in each period the universal constraint $\sum_{k\in[N]} q^t_k =1$ must be fulfilled and thus in each period $t \in [T]$ the universal variables have to obey the domain 
$
\D_All_t =\{\pmb{q}^t \in \{0,1\}^N \mid \sum_{k\in[N]} q_k^t=1\}.$
The cost for edge $(i,j) \in E$ in scenario $k$ and period $t$ is given by $c_{i,j,k}^t \in \mathbb{R}_+$. As before, auxiliary variables $z_t$ are used to bundle the costs incurred in period $t$ and to avoid a nonlinear term in the objective function. The \QIPPU model for this multistage assignment problem is given below.
\begin{subequations}
\begin{align}
\text{\hspace*{-1cm}(\ASSQPU)\hspace*{.5cm}}
\min \  & \sum_{i\in[n]} \sum_{j\in[n]} c_{i,j}^0 x^0_{i,j} +  \sum_{t\in[T]} z_t \label{Eq::Ex_Mod_Ass1}\\
\textnormal{s.t.}\ &\rlap{$\exists \pmb{x}^0 \in \{0,1\}^{n\times n} \quad \forall \pmb{q}^1 \in \D_All_1  \quad \exists \pmb{x}^1 \in \{0,1\}^{n\times n}\cdots$} \nonumber\\
&\rlap{$\cdots  \ \forall \pmb{q}^T \in \D_All_T  \quad \exists \pmb{x}^T \in \{0,1\}^{n\times n}\quad \exists \pmb{z} \in \mathbb{R}_+^T:$}  \nonumber\\
& \sum_{j\in[n]} \sum_{t=0}^T  x_{i,j}^t = 1 && \forall \, i \in [n]\ \label{Eq::Ex_Mod_Ass3}\\
& \sum_{i\in[n]} \sum_{t=0}^T x_{i,j}^t = 1 && \forall \, j \in [n]\ \label{Eq::Ex_Mod_Ass4}\\
& \sum_{i\in[n]}  \sum_{j\in[n]} c_{i,j,k}^t x_{i,j}^t  \leq z_t + M_k^t(1-q_{k}^t) && \hfil \forall \, k \in [N],\,  t \in [T]  \label{Eq::Ex_Mod_Ass5}
\end{align}
\end{subequations}
The Objective~\eqref{Eq::Ex_Mod_Ass1} consists of the expenses from the initial stage with fixed costs and each period with uncertain costs. Constraints~\eqref{Eq::Ex_Mod_Ass3} and  \eqref{Eq::Ex_Mod_Ass4} ensure that the found solution is indeed a perfect matching. Constraint~\eqref{Eq::Ex_Mod_Ass5} linearizes the dependence between selected scenario and incurred costs. Similar to the multistage selection problem, in order to build an equivalent QIP we represent the universal player's decision as an integer variable $\ell_t \in [N]$ and then convert it into existential indicator variables $\pmb{q}^t \in \{0,1\}^N$. We then define the following problem.
\begin{subequations}
\begin{align}
\text{\hspace*{-1cm}(\ASSQ)\hspace*{.5cm}}
\min \  & \sum_{i\in[n]} \sum_{j\in[n]} c_{i,j}^0 x^0_{i,j} +  \sum_{t\in[T]} z_t \label{Eq::Ex_Mod_QAss1}\\
\textnormal{s.t.}\ &\rlap{$\exists \pmb{x}^0 \in \{0,1\}^{n\times n}\quad \forall \ell_1 \in [N] \quad \exists \pmb{q}^1 \in \{0,1\}^N \quad \exists \pmb{x}^1 \in \{0,1\}^{n\times n} \ \cdots  $}~ \nonumber\\
&\rlap{$\cdots \  \forall \ell_T \in [N] \quad \exists \pmb{q}^T \in \{0,1\}^N \quad \exists \pmb{x}^T \in \{0,1\}^{n\times n} \ \exists \pmb{z} \in \mathbb{R}_+^T$:}~ \nonumber\\
& \sum_{j\in[n]} \sum_{t=0}^T  x_{i,j}^t = 1 && \forall \, i \in [n]\ \label{Eq::Ex_Mod_QAss3}\\
& \sum_{i\in[n]} \sum_{t=0}^T x_{i,j}^t = 1 && \forall \, j \in [n]\ \label{Eq::Ex_Mod_QAss4}\\
& \sum_{i\in[n]}  \sum_{j\in[n]} c_{i,j,k}^t x_{i,j}^t  \leq z_t + M_k^t(1-q_{k}^t) && \hfil \forall \, k \in [N],\, t \in [T] \label{Eq::Ex_Mod_QAss5}\\
& \sum_{k\in[N]} q^t_k=1 && \forall \, t\in[T]  \label{Eq::Ex_Mod_QAss6}\\
& \sum_{k\in[N]} k \cdot q^t_k=\ell_t && \forall \, t\in[T]  \label{Eq::Ex_Mod_QAss7}
\end{align}
\end{subequations}

As before, we are interested in a robust counterpart that can be solved using standard MIP solvers. Similar to the notation used in the previous subsection, let $\cR$ denote the set of all possible sequences of scenarios and let $\pmb{r}(t)$ denote the sub-sequence of the scenario sequence $\pmb{r}$ up to period $t$. The variable $x_{i,j}^{\pmb{r}(t)}$ indicates the decision of selecting edge $(i,j)$ in period $t$ after the sub-sequence $\pmb{r}(t)$ of $\pmb{r}$ occurred. The robust counterpart of \ASSQPU is given below.
\begin{subequations}
\begin{align}
\text{\hspace*{-1cm}(\ASSRC)\hspace*{.5cm}}
\min \  & \sum_{i\in[n]} \sum_{j\in[n]} c_{i,j}^0 x^0_{i,j} +  z \label{Eq::Ex_Mod_RCAss1}\\
\textnormal{s.t.}\ &z \geq \sum_{i\in[n]} \sum_{j\in[n]} \sum_{t\in[T]} c_{i,j,r_t}^t x_{i,j}^{\pmb{r}(t)} && \forall\, \pmb{r} \in \cR \label{Eq::Ex_Mod_RCAss2}\\
& \sum_{j\in[n]} \sum_{t=0}^T x_{i,j}^{\pmb{r}(t)} = 1 && \forall \, i \in [n], \,  \pmb{r} \in \cR \label{Eq::Ex_Mod_RCAss3}\\
& \sum_{j\in[n]} \sum_{t=0}^T x_{i,j}^{\pmb{r}(t)} = 1 && \forall \, j \in [n],\,  \pmb{r} \in \cR \label{Eq::Ex_Mod_RCAss4}\\
& x_{i,j}^0 \in\{0,1\} && \forall\, i,\, j \in [n] \label{Eq::Ex_Mod_RCAss5}\\
& x^{\pmb{r}(t)}_{i,j} \in\{0,1\} && \forall \, i,\,j\in [n], \, \pmb{r} \in \cR,\, t \in [T] \label{Eq::Ex_Mod_RCAss6}\\
&z \in\mathbb{R}
\end{align}
\end{subequations}

 \subsection{Multistage Robust Lot-Sizing Problems}\label{SubSec::LotSizing}

We also consider a single item lot-sizing problem with discrete ordering decisions. The model is inspired by the ones used in \cite{bertsimas2015design,bertsimas2018binary}.  After the initial time period there are $T$ more time periods.  At the beginning of each time period $t$, the product demand that needs to be satisfied is disclosed. Here, the lower and upper bounds of the anticipated demand intervals  $[\munderbar{d}_t,\bar{d}_t]$ are considered. Binary universal variables $z_t$ indicate the occurring demand. In order to serve the demand, in each period $t \in \{0,\ldots,T-1\}$ one of $B$ basic orders can be placed resulting in the deliverance of quantity $q_b$ at the beginning of the subsequent period for a unit cost $c_B$. Binary (existential) variables $x^t_{b}$ indicate whether basic order $b$ is selected in period $t$. Additionally, one of $U$ urgent orders of quantity $p_u$ can be placed in each period $t \in \{1,\ldots,T\}$ for a unit cost of $c_N$ with $c_N > c_B$, which is delivered in the same period. If the available quantity exceeds the demand, excess units are stored incurring unit cost $c_S$. Binary (existential) variables $y^{t}_{u}$ indicate whether urgent order $u$ is selected in period $t$. The inventory level in each period is given by $I_t$, with $I_0$ being set to zero. The inventory $I_t$ at period $t \in [T]$ is then given by 
$$I_t=\sum_{t'\leq t} \left( \sum_{b \in [B]} q_b x^{t'-1}_{b}  +\sum_{u \in [U]}  p_u y^{t'}_{u} -(\munderbar{d}_{t'}+\left( \bar{d}_{t'}-\munderbar{d}_{t'})z_{t'}\right) \right)$$
which we explicitly substitute in our implementation.
\begin{subequations}
\begin{align}
\text{\hspace*{-1cm}(\LOTQ)\hspace*{.5cm}}
\min\ & \sum_{t \in [T]} \left( \sum_{b \in [B]} c_B q_b x^{t-1}_{b} + c_S I_t +\sum_{u \in [U]} c_U p_u y^t_{u} \right) \\
\textnormal{s.t. } &  \exists\ \pmb{x}^0 \in \{0,1\}^B \quad \forall\, {z}_1 \in\{0,1\} \quad \exists \, \pmb{y}^1 \in \{0,1\}^U\nonumber\\
&\exists\ \pmb{x}^1 \in \{0,1\}^B  \quad \forall\, {z}_2 \in\{0,1\} \quad \ldots  \nonumber\\
&\exists\ \pmb{x}^{T-1} \in \{0,1\}^B \quad\forall\, {z}_T \in\{0,1\} \quad \exists \, \pmb{y}^T \in \{0,1\}^U: \nonumber\\
&I_t \geq 0\qquad \qquad\forall  t \in [T] 
\end{align}
\end{subequations}

The objective aims at minimizing the overall costs consisting of storage cost as well as basic and urgent order costs, under the constraint that the inventory in each period is non-negative.

Let $\cR=\{0,1\}^T$ denote the set of all possible scenarios and let $\pmb{r}(t)$ denote the sub-sequence of scenario $\pmb{r}$ up to period $t$. The variables $x_{b}^{\pmb{r}(t)}$, $y_{u}^{\pmb{r}(t)}$ and  $I_{\pmb{r}(t)}$ indicate the basic and urgent order decisions if sub-sequence $\pmb{r}(t)$ of $\pmb{r}$ occurred. Let $d_{r_{t}}$ be the demand according to the corresponding scenario. Omitting the binary variable domain for each variable, the corresponding robust counterpart is given below. The variables $I_{\pmb{r}(t)}$ have been substituted in the implementation.
\begin{subequations}
\begin{align}
\text{\hspace*{-.6cm}(\LOTRC)\hspace*{.1cm}}
\min\ & z\\
\textnormal{s.t. } &z\geq \sum_{t \in [T]} \left( \sum_{b \in [B]} c_B q_b x_{b}^{\pmb{r}(t-1)} + c_S I_{\pmb{r}(t)} +\sum_{u\in [U]} c_U p_u y_{u}^{\pmb{r}(t)} \right)\hspace*{-.5cm}~&& \forall \pmb{r} \in \cR \\ 
%\textnormal{s.t. } &z\geq \sum_{\substack{t \in [T]\\b \in [B]}} c_B q_b x_{b}^{\pmb{r}(t-1)} + c_S I_{\pmb{r}(t)} +\sum_{u\in [U]} c_U p_u y_{u}^{\pmb{r}(t)}&& \forall \pmb{r} \in \cR \\ 
&I_{\pmb{r}(t)}=\sum_{t'\leq t} \left( \sum_{b \in [B]} q_b x_{b}^{\pmb{r}(t'-1)}  +\sum_{u \in [U]}  p_u y_{u}^{\pmb{r}(t')} -d_{r_{t}} \right)  &&\forall \pmb{r} \in \cR, \,  t \in [T]\\
&I_{\pmb{r}(t)} \geq 0&&\forall \pmb{r} \in \cR, \,  t \in [T]
\end{align}
\end{subequations}

 \subsection{Multistage Robust Knapsack Problems}\label{SubSec::Knapsack}
 
Finally, we consider a multistage robust knapsack problem. It is based on the multistage model used in~\cite{bampis_et_al:LIPIcs:2019:10966}, where uncertainty in the item weights has been added. 
For each time step $t \in \{0,1,\ldots,T\}$ a knapsack instance with capacity $c_t$, $n$ items, and uncertain weights must be solved. For period $t$ 
we denote by $p^t_{i}$ the profit and by $w^t_{i}$ the basic weight of item $i$. We use $x^t_{i}$ as the binary existential indicator variable whether item $i$ is packed in period $t$. After the initial deterministic knapsack in period $0$, the basic weight $w^t_{i}$ for each $t\in[T]$ is increased by an additional weight $a^t_{i}$ if the corresponding binary universal variable $z^t_{i}$ is set to $1$. Similar to the concept of budgeted uncertainty~\cite{bertsimas2003robust}, in each time step the weight of at most $\alpha$ items can be increased and overall at most $\beta$ such increases are allowed. Hence, the domain of the universal variables depend on the realization of previous universal variables and is given by
$$\mathcal{D}_t\left(\pmb{z}^1,\ldots,\pmb{z}^{t-1}\right)=\left\{\pmb{z}^t \in \{0,1\}^n : \Vert\pmb{z}^t\Vert_1\leq \alpha, \sum_{t'=1}^{t}  \Vert\pmb{z}^{t'}\Vert_1\leq \beta\right\}\, .$$
 Additionally to the profit of each item a transition bonus is used to aim for a stable sequence of solutions: if the selection decision for item $i$ remains the same in period $t$ and $t+1$ a bonus of $b^t_{i}$ is granted, indicated by the existential variable $y^t_{i}$. 
 \begin{subequations}
\begin{align}
\text{\hspace*{-1cm}(\KNAQPU)\hspace*{.5cm}}
\max\ & \sum_{t =0}^T \sum_{i \in [n]} p^t_{i}x^t_{i}+ \sum_{t \in [T]}  \sum_{i \in [n]} b^t_{i}y^t_{i} \\
\textnormal{s.t. } &  
\rlap{$\exists\ \pmb{x}^0 \in \{0,1\}^n \quad \forall\, \pmb{z}^1 \in \mathcal{D}_1\quad   \exists\ \pmb{x}^1 \in \{0,1\}^n  $}~ \nonumber\\
&\rlap{$\exists\ \pmb{y}^1 \in \{0,1\}^n   \quad \forall\, \pmb{z}^2 \in \mathcal{D}_2\quad   \exists\ \pmb{x}^2 \in \{0,1\}^n$}~    \nonumber\\
&\rlap{$\ldots  \forall\, \pmb{z}^T \in \mathcal{D}_T \quad   \exists\ \pmb{x}^T \in \{0,1\}^n\ \exists\ \pmb{y}^T \in \{0,1\}^n   :$}~\nonumber\\
& \sum_{i \in [n]} w^0_{i} x^0_{i} \leq c_0\label{Constraint::Knap_Knapsack1}\\
& \sum_{i \in [n]} \left(w^t_{i} + a^t_{i}z^t_{i} \right)x^t_{i} \leq c_t&&\forall t \in [T]\label{Constraint::Knap_Knapsack2}\\
& y^t_{i} \leq -x^{t-1}_{i}+x^t_{i}+1 && \forall t \in [T], \forall i \in [n] \label{Constraint::Knap_Bonus1}\\
& y^t_{i} \leq x^{t-1}_{i}-x^t_{i}+1 && \forall t \in [T], \forall i \in [n]\label{Constraint::Knap_Bonus2}
\end{align}
\end{subequations}

Constraints \eqref{Constraint::Knap_Knapsack1} and \eqref{Constraint::Knap_Knapsack2} ensure that in each period a valid knapsack solution is found. The nonlinearity in Constraint \eqref{Constraint::Knap_Knapsack2} is linearized, resulting in an additional existential binary variable embedded in the same variable block as the corresponding $x^t_{i}$. Constraints \eqref{Constraint::Knap_Bonus1} and \eqref{Constraint::Knap_Bonus2} link the periods, potentially resulting in the bonus payment. 

Again, let $\cR$ denote the set of all possible scenarios  (realizations of all universal variables) according to the budgeted uncertainty set and let $\pmb{r}(t)$ denote the sub-sequence of scenario $\pmb{r}$ up to period $t$. The variables $x_{i}^{\pmb{r}(t)}$ and $y_{i}^{\pmb{r}(t)}$ indicate the selection and bonus indicators if sub-sequence $\pmb{r}(t)$ of $\pmb{r}$ occurred. Let $w_{i}^{r_{t}}$ be the weight according to the corresponding scenario. The corresponding robust counterpart is given as follows, where we omit stating the binary variable domain for each variable.
\begin{subequations}
\begin{align}
\text{\hspace*{-1cm}(\KNARC)\hspace*{.5cm}}
\max\ & p^0_{i}x^0_{i}+z \\
\textnormal{s.t. } & z \leq \sum_{t \in [T]}\sum_{i \in [n]} p^t_{i}x_{i}^{\pmb{r}(t)} + \sum_{t \in [T]}  \sum_{i \in [n]} b^t_{i}y_{i}^{\pmb{r}(t)}  && \forall \pmb{r} \in \cR\\
& \sum_{i \in [n]} w^0_{i} x^0_{i} \leq c_0\\
& \sum_{i \in [n]}  w_{i}^{r_{t}} x_{i}^{\pmb{r}(t)} \leq c_t&&\forall t \in [T],\pmb{r} \in \cR \\
& y_{i}^{\pmb{r}(t)} \leq -x_{i}^{\pmb{r}(t-1)}+x_{i}^{\pmb{r}(t)}+1 && \forall t \in [T], i \in [n], \pmb{r} \in \cR \\
& y_{i}^{\pmb{r}(t)} \leq x_{i}^{\pmb{r}(t-1)}-x_{i}^{\pmb{r}(t)}+1 && \forall t \in [T], i \in [n], \pmb{r} \in \cR 
\end{align}
\end{subequations}

\section{Computational Experiments}\label{sec:exp}

We conduct experiments on all multistage robust problems presented in Section~\ref{sec:problems}. The aim of these experiments is to answer the questions: To what size can the approaches presented in this paper solve multistage problems to optimality? And which of the exact approaches performs best? 

\subsection{Setup}

We use CPLEX (12.9.0) as MIP solver in order to solve the robust counterpart DEP. The QMIP solver Yasol uses CPLEX (12.6.1) as its black-box LP solver. Since Yasol currently only uses a single thread we also restricted CPLEX to a single thread in order to obtain a balanced comparison (CPLEX turned out to be even faster if restricted to a single thread). All experiments were run on an Intel(R) Core(TM) i7-4790 with 3.60 GHz and 32 GB RAM.

 For each generated instance random values are created using the \texttt{C++} function \texttt{rand()} from the standard general utilities library and the modulo operator. The time limit for each instance is set to $1800$ seconds. Note that for each instance the parameter $T$ representing the number of time periods is equal to the number of universal decision stages.

 For the evaluation of the experiments we present performance profiles \cite{dolan2002benchmarking}. Therefore, we briefly recall this concept: Let $\mathcal{S}$ be the set of considered solvers, $\mathcal{P}$ the set of instances and $t_{p,s}$ the runtime of solver $s$ on instance $p$. 
We assume $t_{p,s}$ is set to infinity (or large enough) if solver $s$ does not solve instance $p$ within the time limit. The percentage of instances for which the performance ratio of solver $s$ is within a factor $\tau\geq 1$ of the best ratio of all solvers is given by
$$p_s(\tau)= \frac{1}{\vert \mathcal{P} \vert} \left\vert \left\{p \in \mathcal{P} \mid\ \frac{t_{p,s}}{\min\limits_{\hat{s}\in \mathcal{S}}t_{p,\hat{s}}} \leq \tau \right\} \right\vert \, .$$
Hence, the function $p_s$ can be viewed as the distribution function for the performance ratio, which is plotted in a performance profile for each solver. For each  performance profile in this paper we calculate $p_s(\tau)$ for $\tau =1+0.5t$ with $t \in \mathbb{N}_0$ and $t$ large enough until $p_s(\tau)$ is constant.  Note that for each instance with (rounded) runtime of $0$ seconds we used the runtime of $1$ second in order to be able to generate useful performance profiles. In each performance profile we use the colors red, blue and black for the DEP, QIP and \QIPPU, respectively.

\subsection{Experiments on Multistage Robust Selection}\label{sec:expsel}

We start with the multistage selection problem as introduced in Section \ref{SubSec::Selection} in which $p$ out of $n$ items  must be selected. We compare the performance of Yasol on the quantified models \SELQPU and \SELQ with the performance of CPLEX on the robust counterpart \SELRC.
An instance is given by the number of available items $n$, the number of items $p$ to be selected, the number of time periods $T$ and the number of scenarios $N$ per period.
Recall that $T$ denotes the number of $\max\min$ blocks after the first $\min$. Hence, $T=1$ corresponds to two-stage robust optimization.
We limit experiments to instances with $n=2p$ and thus, the value of $p$ is omitted from now on. The remaining parameters of an instance are the costs $c_{i,k}^t$ of each each item $i$ in scenario $k$ of period $t$, which are randomly selected from the range \{$0,1,\ldots,99\}$. 
 
We compare the performance of solution methods for an increasing number of scenarios.
To this end, we fix the  number of items to $n=10$. This might seem small but it necessary  in order to allow large values of $T$ and $N$ and still be able to find the optimal solution for many instances in reasonable time. 
For most constellations of $T \in \{1,\ldots,8\}$ and $N\in \{2^1,\ldots,2^8\}$, $50$ instances are created. Not all constellations are considered, since the prospects of finding the optimal solution within the time limit of $1800$ seconds is small if both $T$ and $N$ are large. In Appendix~\ref{sec:comp-N} we provide additional results on experiments with fixed number of scenarios $N=4$ and varying  numbers of items $n\in\{10,20,30,40,50\}$ and periods $T\in \{1,\ldots,8\}$, where findings are similar to what is presented here.

Table~\ref{Table::NumberSolvedSelectionTriange} shows the number of instances solved. We omit the row for $N=2$ and the column for $T=1$, as solving those instance for each model never took longer than two seconds. For cells marked with a hyphen no experiments were conducted, as it is expected that none (or very few) instances would be solved in the given time limit. The zeros in brackets indicate that the robust counterparts could not be created due to their size exceeding 200GB. 
\begin{table}[ht]
\footnotesize
\centering
\caption[Number of solved selection instances with $n=10$]{Number of solved selection instances with $n=10$ items and various $N$ and $T$.}\label{Table::NumberSolvedSelectionTriange}
\begin{tabular}{lL{1.4cm}R{1.1cm}R{1.1cm}R{1.1cm}R{1.1cm}R{1.1cm}R{1.1cm}R{1.1cm}R{1.1cm}}

\toprule
&model & $T=2$ & $T=3$ & $T=4$ & $T=5$ & $T=6$ & $T=7$ & $T=8$      \\
 \midrule  \multirow{3}*{$N=2^2$}    & \SELRC      & 50     & 50     & 50     & 50     & 50     & 49     & 26 \\
       & \SELQ       & 50     & 50     & 50     & 50     & 50     & 50     & 50 \\
       & \SELQPU      & 50     & 50     & 50     & 50     & 50     & 50     & 50 \\
 \midrule  \multirow{3}*{$N=2^3$}   & \SELRC       & 50     & 50     & 50     & 48     & 4      & 0      & (0) \\
       & \SELQ       & 50     & 50     & 50     & 50     & 50     & 47     & 0  \\
       & \SELQPU       & 50     & 50     & 50     & 50     & 50     & 50     & 38 \\
 \midrule  \multirow{3}*{$N=2^4$}   & \SELRC      & 50     & 50     & 44     & 0      & (0)     & (0)     & -  \\
       &\SELQ      & 50     & 50     & 50     & 50     & 0      & 0      & -  \\
       & \SELQPU      & 50     & 50     & 50     & 50     & 36     & 8      & -  \\
 \midrule  \multirow{3}*{$N=2^5$}    & \SELRC      & 50     & 50     & 3      & (0)     & (0)     & -      & -  \\
       &\SELQ     & 50     & 50     & 50     & 0      & 0      & -      & -  \\
       & \SELQPU      & 50     & 50     & 50     & 23     & 1      & -      & -  \\
  \midrule \multirow{3}*{$N=2^6$}   & \SELRC     & 50     & 32     & (0)     & (0)     & -      & -      & -  \\
       &\SELQ      & 50     & 50     & 0      & 0      & -      & -      & -  \\
       & \SELQPU       & 50     & 50     & 27     & 2      & -      & -      & -  \\
 \midrule  \multirow{3}*{$N=2^7$}  & \SELRC      & 50     & 0      & (0)     & -      & -      & -      & -  \\
       &\SELQ     & 50     & 14     & 0      & -      & -      & -      & -  \\
       & \SELQPU      & 50     & 48     & 3      & -      & -      & -      & -  \\
 \midrule  \multirow{3}*{$N=2^8$}   & \SELRC       & 48     & (0)     & -      & -      & -      & -      & -  \\
       &\SELQ      & 50     & 0      & -      & -      & -      & -      & -  \\
       & \SELQPU      & 50     & 11     & -      & -      & -      & -      & - \\\bottomrule
\end{tabular}
\end{table}

For increasing $N$ and $T$ the number of quantified programs solved by Yasol tends to be larger than the number of robust counterpart solved by CPLEX. For each configuration, the number of solved  \SELQPU models  is highest and at least one \SELQPU instance is always solved. On $243$ of the $303$ instances where no model was solved to optimality, \SELQPU resulted in the best incumbent solution. On $175$ of those instances \SELQPU was the only model for which any solution was found at all. On additional $23$ instances the optimal solution of \SELQPU was found while for the other models not even an incumbent solution was found. 
In Table \ref{Table::TimeSelectionTriange} the average runtimes are presented. Bold values represent cases where each instance reached the time limit.

\begin{table}[ht]
\footnotesize
\centering
\caption{Average runtime of selection instances with $n=10$.}\label{Table::TimeSelectionTriange}
\begin{tabular}{llR{1.1cm}R{1.1cm}R{1.1cm}R{1.1cm}R{1.1cm}R{1.1cm}R{1.1cm}R{1.1cm}}

\toprule
&model    & $T=2$ & $T=3$ & $T=4$ & $T=5$ & $T=6$ & $T=7$ & $T=8$      \\
  \midrule \multirow{3}*{$N=2^2$}  & \SELRC  &0.0       & 0.0         & 0.1    & 1.5     & 10.8    & 223.3 & 1396.3 \\
      &\SELQ   &0.1   & 0.2       & 0.2    & 0.7    & 2.1    & 8.6   & 33.1     \\
& \SELQPU    &0.1     & 0.2      & 0.2    & 0.5    & 1.1    & 4.1   & 11.9     \\
  \midrule \multirow{3}*{$N=2^3$}  & \SELRC&  0.0       & 0.3       & 7.3    & 366.0  & 1756.9 & \textbf{1800.0}   & \textbf{1800.0}       \\
      &\SELQ   &0.2    & 0.4       & 2.0    & 14.0   & 103.86  & 905.3  & \textbf{1800.0}      \\
& \SELQPU    &0.1   & 0.4      & 1.0       & 4.0      & 26.2   & 146.8 & 947.9       \\
  \midrule \multirow{3}*{$N=2^4$}  & \SELRC & 0.0       & 6.3      & 564.7  & \textbf{1800.0}    & \textbf{1800.0}    & \textbf{1800.0}   & -           \\
      &\SELQ   &0.2     & 2.2      & 30.1   & 557.5  & \textbf{1800.0}    & \textbf{1800.0}   & -             \\
& \SELQPU    &0.3   & 1.1     & 7.9    & 111.7  & 1092.2  & 1655.2 & -            \\
  \midrule \multirow{3}*{$N=2^5$}  & \SELRC & 0.8    & 102.4    & 1792.1 & \textbf{1800.0}    & \textbf{1800.0}    & -      &       \\
      &\SELQ   &0.8    & 19.7     & 600.0     & \textbf{1800.0}    & \textbf{1800.0}    & -      & -            \\
& \SELQPU    &0.5    & 5.3      & 129.3   & 1268.5  & 1783.4  & -      & -            \\
  \midrule \multirow{3}*{$N=2^6$}  & \SELRC & 4.1    & 1164.8    & \textbf{1800.0}    & \textbf{1800.0}    & -       & -      & -            \\
      &\SELQ   &3.6    & 237.3     & \textbf{1800.0}    & \textbf{1800.0}    & -       & -      & -           \\
& \SELQPU    &1.4   & 51.4     & 1260.5 & 1787.6 & -       & -      & -           \\
  \midrule \multirow{3}*{$N=2^7$}  & \SELRC & 25.4   & 1829.1 & \textbf{1800.0}    & -       & -       & -      & -           \\
      &\SELQ   &20.0   & 1638.1   & \textbf{1800.0}    & -       & -       & -      & -           \\
& \SELQPU    &6.8    & 415.9    & 1715.9 & -       & -       & -      & -           \\
  \midrule \multirow{3}*{$N=2^8$}  & \SELRC & 252.1 & \textbf{1800.0}         & -       & -       & -       & -      & -           \\
      &\SELQ   &154.5  & \textbf{1800.0}      & -       & -       & -       & -      & -           \\
& \SELQPU    &51.2   & 1540.6   & -       & -       & -       & -      & -           \\\bottomrule
\end{tabular}
\end{table}
We highlight the modest increase of the runtime of \SELQPU models for increasing $T$ and $N$, even compared to \SELQ. Tables~\ref{Table::NumberSolvedSelectionTriange} and~\ref{Table::TimeSelectionTriange} show that a) for instances with many periods and scenarios the use of quantified programs is superior to solving the robust counterpart, and b) utilizing universal constraints rather than standard QIPs  is of  advantage in this setting. 

\begin{figure}[h!]
\centering
\includegraphics[scale=1]{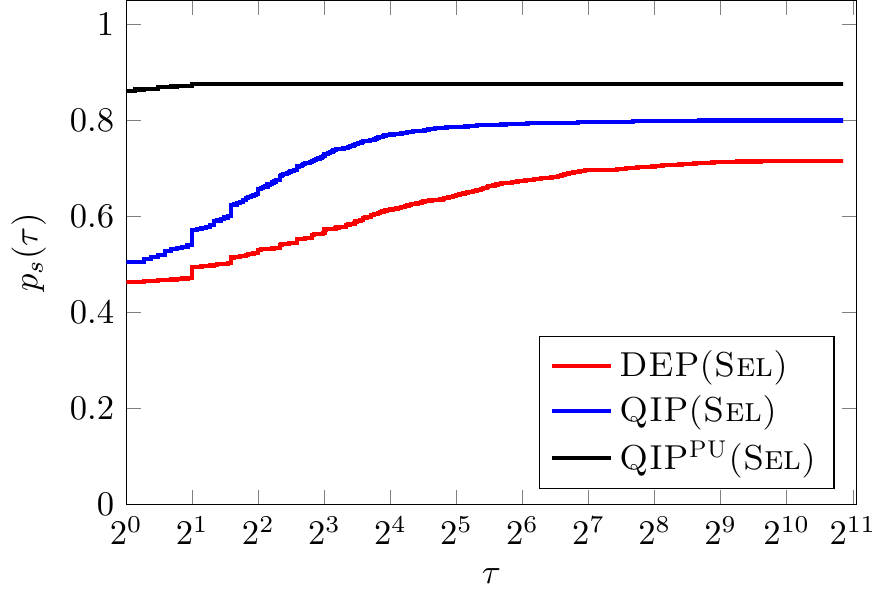}
\caption{Performance profile for all selection instances with $n=10$.}\label{Fig::PerformanceSelectionAlln10}
\end{figure}

The performance profile in Figure \ref{Fig::PerformanceSelectionAlln10} summarizes these results. It shows that no other model is solved faster than then \SELQPU model in $86\%$ of the instances. $16\%$ more instances were solved when modeled as \SELQPU  than \SELRC.

\subsection{Experiments on Multistage Robust Assignment}\label{sec:expass}

The multistage assignment problem as introduced in Subsection \ref{SubSec::Assignment} deals with the problem of finding a perfect matching in a bipartite graph with minimal costs. We compare the performance of Yasol on the quantified models \ASSQPU and \ASSQ with the performance of CPLEX on the robust counterpart \ASSRC.

Each instance is given by the size $n$ of each partition, the number of time periods $T$ and the number of scenarios $N$ per period. The remaining parameters of an instance are the cost $c_{i,j,k}^t$ for each edge $(i,j)$ in scenario $k$ of period $t$, which are randomly selected from the range $\{0,1,\ldots,99\}$. 

For each $n \in\{4,\ldots,10\}$, $T \in \{1,\ldots,4\}$ and $N\in \{2,4,8\}$ we generated $50$ instances for each model type \ASSQPU, \ASSQ and \ASSRC. As before, each solver has a time limit of $1800$ seconds. We examine how the realizations of $n$, $T$ and $N$ affects the runtime and which model-solver combination is best suited for the different instances. 

In Table \ref{Table::AssignmentNumTimeN4} the number of solved instances as well as the overall average runtime is presented for instances with $N=4$ scenarios. 
For instances with $N=2$ and $N=8$ scenarios the respective Tables~\ref{Table::AssignmentNumTimeN2} and~\ref{Table::AssignmentNumTimeN8} can be found in Appendix~\ref{Appendix::Assignment}.

\begin{table}[ht]
\centering
\caption{Number of solved assignment instances (opt) with $N=4$ and the average runtime (time).}\label{Table::AssignmentNumTimeN4}
\footnotesize
\begin{tabular}{lp{1cm}R{.9cm}rR{.9cm}rR{.9cm}rR{.9cm}r}
\toprule
& & \multicolumn{2}{c}{\hspace{.45cm}$T=1$}    &\multicolumn{2}{c}{\hspace{.45cm}$T=2$}&\multicolumn{2}{c}{\hspace{.45cm}$T=3$}&\multicolumn{2}{c}{\hspace{.45cm}$T=4$}    \\
& model& opt &time  & opt &time & opt &time & opt &time     \\
 \midrule \multirow{3}*{$n=4$}  & \ASSRC  & 50     &$ 0.0 $ & 50     &$ 0.0   $ & 50     &$ 0.0   $ & 50     &$ 3.1    $ \\
     & \ASSQ                & 50     &$ 0.1 $ & 50     &$ 0.1   $ & 50     &$ 0.3   $ & 50     &$ 0.7    $ \\
     & \ASSQPU & 50     &$ 0.1 $ & 50     &$ 0.1   $ & 50     &$ 0.2   $ & 50     &$ 0.5    $ \\
\midrule \multirow{3}*{$n=5$}    & \ASSRC   & 50     &$ 0.0 $ & 50     &$ 0.0   $ & 50     &$ 0.6   $ & 50     &$ 14.5   $ \\
     & \ASSQ                 & 50     &$ 0.0 $ & 50     &$ 0.3   $ & 50     &$ 1.0   $ & 50     &$ 2.6    $ \\
     & \ASSQPU & 50     &$ 0.1 $ & 50     &$ 0.2   $ & 50     &$ 0.5   $ & 50     &$ 1.8    $ \\
\midrule \multirow{3}*{$n=6$}    & \ASSRC  & 50     &$ 0.0 $ & 50     &$ 0.0   $ & 50     &$ 2.0   $ & 50     &$ 46.9   $ \\
     & \ASSQ                 & 50     &$ 0.2 $ & 50     &$ 1.0   $ & 50     &$ 3.9   $ & 50     &$ 13.3   $ \\
     & \ASSQPU & 50     &$ 0.2 $ & 50     &$ 0.6   $ & 50     &$ 1.9   $ & 50     &$ 6.4    $ \\
\midrule \multirow{3}*{$n=7$}    & \ASSRC  & 50     &$ 0.0 $ & 50     &$ 0.0   $ & 50     &$ 8.9   $ & 48     &$ 218.6  $ \\
     & \ASSQ                 & 50     &$ 0.2 $ & 50     &$ 5.3   $ & 50     &$ 19.4  $ & 50     &$ 72.7   $ \\
     & \ASSQPU & 50     &$ 0.2 $ & 50     &$ 1.7   $ & 50     &$ 6.0   $ & 50     &$ 25.2   $ \\
\midrule \multirow{3}*{$n=8$}   & \ASSRC  & 50     &$ 0.0 $ & 50     &$ 0.1   $ & 50     &$ 26.1  $ & 33     &$ 922.7  $ \\
     & \ASSQ                 & 50     &$ 0.4 $ & 50     &$ 14.4  $ & 50     &$ 76.3  $ & 50     &$ 257.5  $ \\
     & \ASSQPU & 50     &$ 0.3 $ & 50     &$ 3.5   $ & 50     &$ 26.0  $ & 50     &$ 117.0  $ \\
\midrule \multirow{3}*{$n=9$}   & \ASSRC  & 50     &$ 0.0 $ & 50     &$ 0.7   $ & 50     &$ 155.5 $ & 15     &$ 1444.8 $ \\
     & \ASSQ                 & 50     &$ 0.6 $ & 50     &$ 59.9  $ & 50     &$ 365.2 $ & 33     &$ 1318.8 $ \\
     & \ASSQPU & 50     &$ 0.4 $ & 50     &$ 11.2  $ & 50     &$ 102.9 $ & 49     &$ 554.2$ \\
\midrule \multirow{3}*{$n=10$}   & \ASSRC & 50     &$ 0.0 $ & 50     &$ 1.0   $ & 47     &$ 318.8 $ & 8      &$ 1658.0  $ \\
     & \ASSQ                 & 50     &$ 1.4 $ & 50     &$ 192.9 $ & 36     &$ 1217.5 $ & 1      &$ 1792.8 $ \\
     & \ASSQPU & 50     &$ 0.8 $ & 50     &$ 43.4  $ & 48     &$ 580.7 $ & 12     &$ 1689.4 $\\\bottomrule
\end{tabular}
\end{table}

\begin{figure}[h!]
\centering
\includegraphics[scale=1]{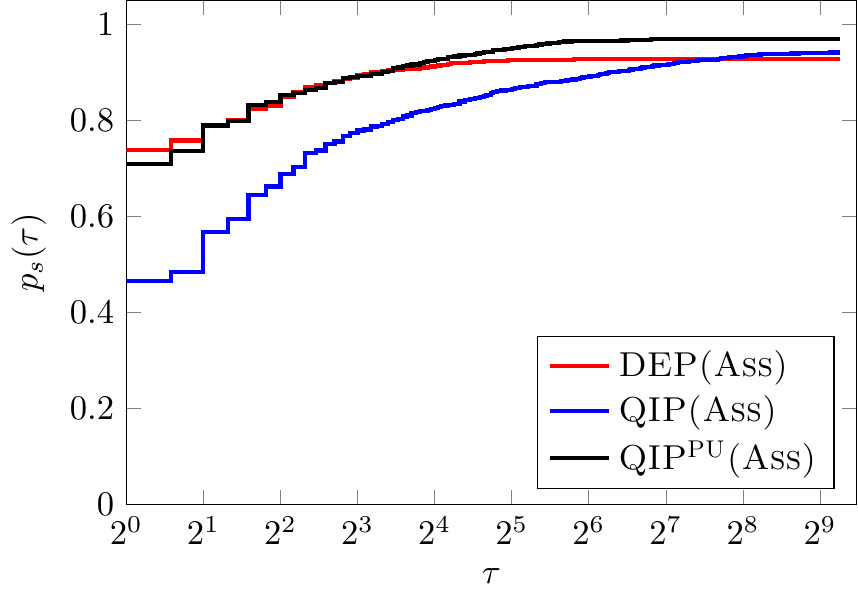}
\caption{Performance profile for all assignment instances with $N=4$.}\label{Fig::PerformanceAssAll}
\end{figure}

For increasing $T$ and $N$, the \ASSQPU model tends to yield the best results but CPLEX (on \ASSRC) remains highly competitive and is able to catch up with increasing $n$. 
In Figure \ref{Fig::PerformanceAssAll} the performance profile on the $1400$ instances with $N=4$ is given.
It can be seen that CPLEX is significantly faster on most \ASSRC instances and is fastest on more than $75\%$ of the instances. However, overall more of the \ASSQPU and \ASSQ instances are solved. 

\subsection{Experiments on Multistage Robust Lot-Sizing}

We conducted experiments for $B \in \{3,4\}$ basic orders, $U \in \{2,3\}$ urgent orders, and a time horizon of $T \in \{5,\ldots,10\}$ with $50$ instances per constellation resulting in $1200$ instances. For our instances we aimed to stay as close as possible to the data used in \cite{bertsimas2015design} while we did not use a cumulative ordering budget.
Basic and urgent ordering costs are chosen from $c_B \in \{0,1,\ldots,5\}$ and $c_U\in \{0,1,\ldots,10\}$ with $c_B<c_U$. Storage costs are chosen from $c_S \in \{0,1,\ldots,10\}$. For the quantities of the basic orders we use $q_1=64$, $q_2=32$, $q_3=16$, and (in case $B=4$) $q_4=8$. For the quantities of the urgent orders we use $p_u=\frac{100}{u}$ for $u \in \{1,\ldots,U\}$.  The lower and upper bounds for the demand are chosen from $\munderbar{d}_t \in \{0,\ldots,25\}$ and $\bar{d}_t  \in \{75,\ldots,100\}$. The initial inventory $I_0$ is assumed to be zero.

Table~\ref{Table::NumberSolvedLotSizingAllB} shows the number of instances solved for each parameter choice within the 1800 seconds time limit. In both $B$ and $U$, it can be observed that an increasing value leads to less instances being solved. With increasing $T$, the number of solved instances drops, where Yasol is less affected than CPLEX for the smaller instances with $B=3$ (e.g., 17 instances are solve by \LOTRC for $T=10$ and $U=2$, while \LOTQ solved 29 instances), but this is the other way around for the larger instances with $B=4$, resulting in mixed overall results.

\begin{table}[h!]
\centering
\footnotesize
\caption{Number of solved lot-sizing instances.}\label{Table::NumberSolvedLotSizingAllB}
\begin{tabular}{lllrrrrrr}
 \toprule
&& Solver     & $T=5$ & $T=6$ & $T=7$ & $T=8$ & $T=9$ & $T=10$ \\\midrule
\multirow{2}*{$B=3$} & \multirow{2}*{$U=2$} & \LOTRC & 50  & 50  & 49  & 41  & 30  & 17   \\
 &   & \LOTQ      & 50  & 50  & 50  & 50  & 42  & 29   \\\midrule
\multirow{2}*{$B=3$} & \multirow{2}*{$U=3$} & \LOTRC & 50  & 50  & 49  & 37  & 19  & 9    \\
   & & \LOTQ      & 50  & 50  & 50  & 49  & 29  & 17   \\\midrule\midrule
\multirow{2}*{$B=4$}&\multirow{2}*{$U=2$} & \LOTRC & 50  & 50  & 47  & 33  & 30  & 8    \\
 &   & \LOTQ      & 50  & 50  & 50  & 32  & 14  & 2    \\\midrule
\multirow{2}*{$B=4$}&\multirow{2}*{$U=3$} & \LOTRC & 50  & 50  & 44  & 23  & 13  & 7    \\
&    & \LOTQ      & 50  & 50  & 47  & 26  & 11  & 5   \\\bottomrule
\end{tabular}
\end{table}

The number of solved instances is complemented by Table~\ref{Table::AverageRuntimeLotSizingAllB}, where we show average runtimes. We find a similar trend with \LOTQ performing better for $B=3$ and \LOTRC performing better for $B=4$.

\begin{table}[h!]
\centering
\footnotesize
\caption{Average runtime of lot-sizing instances.}\label{Table::AverageRuntimeLotSizingAllB}
\begin{tabular}{lllrrrrrr}
 \toprule
&& Solver     & $T=5$ & $T=6$ & $T=7$ & $T=8$ & $T=9$ & $T=10$ \\\midrule
\multirow{2}*{$B=3$} & \multirow{2}*{$U=2$} & \LOTRC &0.0 & 0.9  & 58.5  & 423.4  & 858.5  & 1318.3 \\
 &   & \LOTQ & 1.5 & 5.9  & 27.3  & 104.6  & 658.4  & 1111.9 \\\midrule
\multirow{2}*{$B=3$} & \multirow{2}*{$U=3$} & \LOTRC &0.3 & 3.9  & 140.5 & 592.9  & 1240.4 & 1525.0 \\
  & & \LOTQ      & 2.8 & 15.9 & 81.2  & 282.5  & 1121.5 & 1378.5 \\\midrule\midrule
\multirow{2}*{$B=4$}&\multirow{2}*{$U=2$}& \LOTRC &0.1 & 3.3  & 144.4 & 693.3  & 1014.9 & 1517.0 \\
&   & \LOTQ      & 4.4 & 29.8 & 242.9 & 1189.4 & 1537.6 & 1761.7 \\\midrule
\multirow{2}*{$B=4$}&\multirow{2}*{$U=3$}& \LOTRC &0.5 & 8.9  & 370.2 & 1091.0 & 1429.0 & 1549.7 \\
&   & \LOTQ      & 7.4 & 59.6 & 517.6 & 1281.1 & 1518.0 & 1633.0\\\bottomrule
\end{tabular}
\end{table}

\begin{figure}[h!]
\centering
\includegraphics[scale=1]{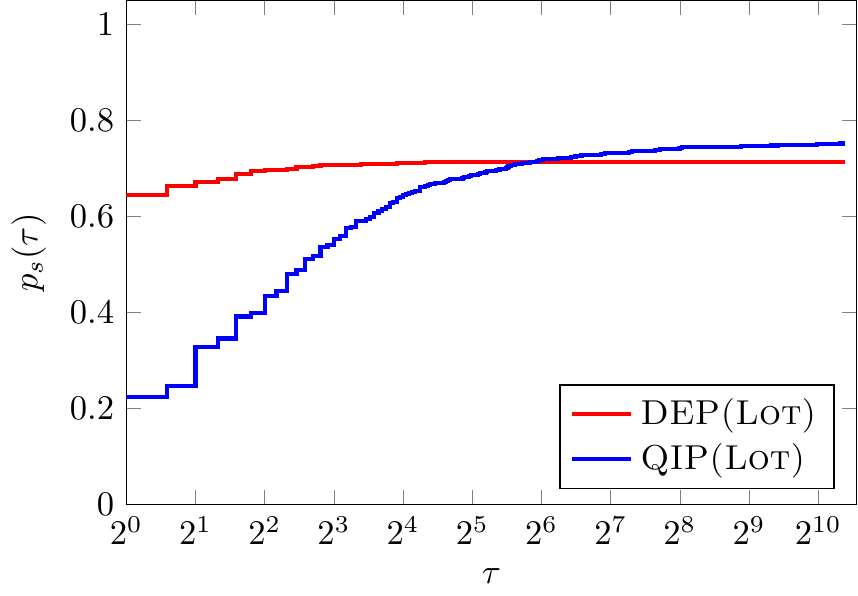}
\caption{Performance profile for all lot-sizing instances.}\label{Fig::PerformanceLotAll}
\end{figure}

Finally, we show a performance profile in Figure~\ref{Fig::PerformanceLotAll}. It can be seen that \LOTRC is the fastest method on more instances than \LOTQ (left-hand side), but \LOTQ is able to solve more instances overall (right-hand side).

\subsection{Experiments on Multistage Robust Knapsack}

In our final experiment, we consider $n \in \{2,\ldots,7\}$ items and a time horizon of $T \in \{1,\ldots,6\}$ with $50$ instances per constellation resulting in $1800$ knapsack instances. We set $\alpha=\lceil{\frac{n}{T+1}}\rceil+1$ and $\beta=n$ which limit the number of items with increased weight per period and over all periods, respectively. For each item $i$ and period $t$ the profits $p^t_{i}$, basic weights $w^t_{i}$, and additional weights $a^t_{i}$ are randomly selected from the ranges $\{0,1,\ldots,100\}$, $\{0,1,\ldots,50\}$, and $\{5,\ldots,20\}$, respectively. The bonuses $b^t_{i}$ are randomly selected from the range $\{0,1,\ldots,50\}$ for each item $i$ and are constant over all periods. The knapsack capacity $c_t$ is also constant over all periods and selected from $\{\lfloor\frac{1}{3}\sum_{i\in [n]}w^0_{i}\rfloor,\ldots,\sum_{i\in [n]}w^0_{i}\}$. 

In Table~\ref{Table::NumberSolvedKnapsack} the number of solved instances for each parameter choice is displayed, where we omit the cases $n \in \{2,\ldots,4\}$, as all instances were solved by both methods. We find that for $n=5$ items, there is a sharp drop in performance for \KNARC when reaching $T=5$, while for $n=6$ and $n=7$, this is already the case for $T=4$ and $T=3$, respectively. When solving instances with \KNAQPU, instances with one or two additional stages can still be solved.

\begin{table}[h!]
\footnotesize
\centering
\caption{Number of solved knapsack instances.}\label{Table::NumberSolvedKnapsack}
\begin{tabular}{llR{1.1cm}R{1.1cm}R{1.1cm}R{1.1cm}R{1.1cm}R{1.1cm}R{1.1cm}R{1.1cm}}

\toprule
&model &$T=1$      & $T=2$ & $T=3$ & $T=4$ & $T=5$ & $T=6$  \\
  \midrule \multirow{2}*{$n=5$} & \KNARC  & 50     & 50     & 50     & 41     & 5     & 3    \\
       & \KNAQPU   & 50     & 50     & 50     & 50     & 48     & 38     \\
\midrule \multirow{2}*{$n=6$} & \KNARC  & 50     & 50     & 25     & 1     & 0     & 0      \\
       & \KNAQPU   & 50     & 50     & 50     & 47     & 18     & 0     \\
\midrule \multirow{2}*{$n=7$} & \KNARC &50 & 47     & 1     & 0     & 0     & 0          \\
       & \KNAQPU   & 50     & 50     & 50     & 5     & 0     & 0     \\\bottomrule
\end{tabular}
\end{table}

\begin{figure}[h!]
\centering
\includegraphics[scale=1]{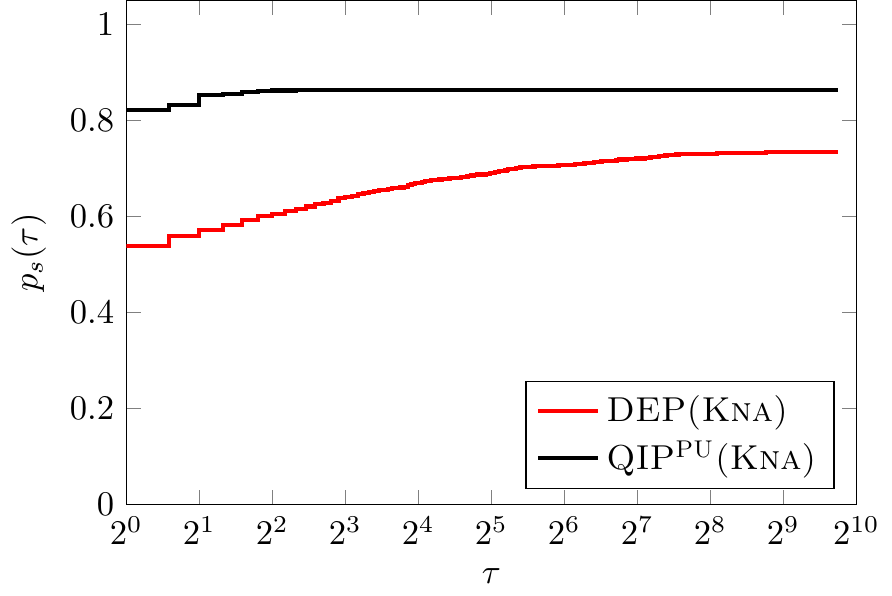}
\caption{Performance profile for all knapsack instances.}\label{Fig::PerformanceKnapAll}
\end{figure}

In Table~\ref{Table::TimeKnapsack} we show the average runtime for each of parameter choice. In all but the very smallest instances, \KNAQPU clearly outperforms \KNARC (up to a factor 100 for $n=7$ and $T=2$). This difference in performance is also recognizable in Figure~\ref{Fig::PerformanceKnapAll}, where the performance profile gives a curve for \KNAQPU that dominates the curve of \KNARC.

\begin{table}[h!]
\footnotesize
\centering
\caption{Average runtime of knapsack instances.}\label{Table::TimeKnapsack}
\begin{tabular}{llR{1.1cm}R{1.1cm}R{1.1cm}R{1.1cm}R{1.1cm}R{1.1cm}R{1.1cm}R{1.1cm}}

\toprule
&model &$T=1$      & $T=2$ & $T=3$ & $T=4$ & $T=5$ & $T=6$  \\
 \midrule \multirow{2}*{$n=2$} & \KNARC  &0.0    &0.0      & 0.0       & 0.0      & 0.0       & 0.0      \\
      & \KNAQPU   &0.0 & 0.0   & 0.1    & 0.1    & 0.14    & 0.3   \\
 \midrule \multirow{2}*{$n=3$} & \KNARC  &0.0    & 0.0     & 0.0       & 0.0   & 0.1    & 1.0   \\
      & \KNAQPU   &0.0 & 0.1   & 0.1    & 0.8   & 1.1    & 2.6   \\
 \midrule \multirow{2}*{$n=4$} & \KNARC  &0.0    & 0.0      & 0.5    & 7.5    & 60.1   & 221.8 \\
      & \KNAQPU   &0.1 & 0.2   & 0.8     & 2.9    & 10.5    & 48.6  \\
 \midrule \multirow{2}*{$n=5$} & \KNARC  &0.0    & 0.7    & 59.8   & 675.7 & 1664.4  & 1712.5 \\
      & \KNAQPU   &0.1 & 0.4   & 3.8    & 26.1  & 277.0  & 874.7 \\
 \midrule \multirow{2}*{$n=6$} & \KNARC  &0.0    & 10.3  & 1344.1 & 1772.3 & \textbf{1800.0}    & \textbf{1800.0}   \\
      & \KNAQPU   &0.1 & 1.3   & 27.5   & 542.5 & 1498.9 & \textbf{1800.0}   \\
 \midrule \multirow{2}*{$n=7$} & \KNARC  &0.0 & 536.0 & 1766.1  & \textbf{1800.0}   & \textbf{1800.0}    & \textbf{1800.0}   \\
      & \KNAQPU   &0.1 & 5.1   & 323.4  & 1702.1 & \textbf{1800.0}    & \textbf{1800.0}   \\\bottomrule
\end{tabular}
\end{table}

\section{Conclusion}\label{sec:conclusion}
%\michael{Nochmal hervorheben, dass wir tatsächlich optimal lösen und keine Decision Rules oder Heuristiken anwenden. Vielleicht auch schon weiter vorne irgendwo; in der Introduction?}
Solving multistage robust problems is a formidable challenge, as the scenario tree grows exponentially. At the same time, from a practical perspective, there is little reason to only consider problems with at most two stages. Hence, new methods to solve multistage robust problems are required.

In this paper we argue that multistage robust discrete  problems can be naturally phrased as quantified integer programs. This allows us to use tools developed for QIPs, in particular the solver Yasol, to solve robust problems to optimality. Two variants of QIP models were presented, which differ regarding the use of a constraint system for the universal player, \IE the modelling capabilities for the uncertainty set.

Our experiments consider a wide range of problems with uncertainty in the objective (selection and assignment), in the right-hand side (lot-sizing), or in the constraints (knapsack).
We find that while CPLEX as an MIP solver for the equivalent deterministic problem performs better for few stages, Yasol has a clear advantage as the number of stages or the number of scenarios per stage increases. In some cases we were not even able to build the deterministic equivalent problem, due to its size. Additionally, the quantified models explicitly containing constraints on universal variables are solved faster than the corresponding QIP without universal constraint system.

These observations have an intuitive explanation. To build the deterministic equivalent problem, we have to take every possible adversarial action into account. On the one hand, for many stages, the size of problems grows exponentially, which means that every additional stage adds a greater burden on CPLEX. In contrast, QIPs can scale better, as the adversarial decisions are considered only as needed. On the other hand, if there are only few stages, CPLEX can be more powerful, as it considers the full list of adversarial actions in one go, whereas Yasol may spend time in the wrong branches of the decision tree while exploring the problem.

There are several interesting further research questions that arise from this work. So far, we have only considered random data. How does the QIP approach perform when using real-world data, e.g., when deriving scenarios using approaches from data-driven robust optimization? Furthermore, is it possible to develop techniques that are specific to QIPs with the structure considered here, as opposed to the more generic capabilities of Yasol? Our encouraging results show that such an approach may lead to an additional performance boost in the solution of multistage robust discrete problems.

\section*{Acknowledgements}\label{sec:Acknowledge}
This research was partially funded by the Deutsche Forschungsgemeinschaft (DFG, German Research Foundation) - 399489083.

\appendix

\section{Additional Results on Multistage Robust Selection}\label{sec:comp-N}

We briefly present results for the selection problem with $N=4$ scenarios in each period and a varying number of items $n\in\{10,20,30,40,50\}$ and number of periods $T\in \{1,\ldots,8\}$. For each setting 50 instances are created and solved with of each of the three models \SELQPU, \SELQ and \SELRC.
For $T<6$ all instances of each model are solved within the time limit of 1800 seconds. In Table \ref{Table::NumberSolvedSelection} the numbers of solved instances for $T\geq 6$ are displayed.
 \begin{table}[h!]
 \centering
 \footnotesize
 \caption{Number of solved selection instances with $N=4$ and $T\in\{6,7,8\}$.}\label{Table::NumberSolvedSelection}
\begin{tabular}{llR{1.1cm}R{1.1cm}R{1.1cm}R{1.1cm}R{1.1cm}R{1.1cm}R{1.1cm}R{1.1cm}}
\toprule
&model & $n=10$    &  $n=20$    &  $n=30$    &  $n=40$    &  $n=50$  \\\midrule
  \multirow{3}*{$T=6$}   & \SELRC  & 50     & 50     & 48     & 50     & 49   \\
       & \SELQ              & 50     & 50     & 50   & 50     & 46    \\
       & \SELQPU & 50     & 50     & 50   & 49     & 48     \\
 \midrule \multirow{3}*{$T=7$}  & \SELRC    & 49     & 32     & 21     &15     & 7   \\
       &  \SELQ                & 50     & 50     & 50     & 49     & 31   \\
       &  \SELQPU & 50     & 50     & 50     & 50     & 33    \\
 \midrule\multirow{3}*{$T=8$}  & \SELRC    & 30     & 0     &0     & 0     & 0      \\
       & \SELQ               & 50     & 50     & 49     & 34     & 15   \\
       &  \SELQPU & 50     & 50     & 50      & 41     & 20   \\\bottomrule
\end{tabular}
\end{table}
In Table \ref{Table::TimeSolvedSelection} the overall average runtime is presented.
\begin{table}[h!]
 \footnotesize
\centering
 \caption{Average runtime of selection instances with $N=4$.} \label{Table::TimeSolvedSelection}
\begin{tabular}{lL{1.4cm}R{1.05cm}R{1.05cm}R{1.05cm}R{1.05cm}R{1.05cm}R{1.05cm}R{1.05cm}R{1.05cm}}
\toprule
&model & $T=1$   & $T=2$ & $T=3$ & $T=4$  & $T=5$  & $T=6$   & $T=7$   & $T=8$            \\
\midrule\multirow{3}*{$n=10$}   & \SELRC &0.0 & 0.0 & 0.0  & 0.0  & 1.5   & 13.3  & 218.9  & 1330.0 \\
  & \SELQ &0.1 & 0.1 & 0.2  & 0.3  & 0.7   & 2.6   & 7.5    & 35.4   \\
     & \SELQPU & 0.1 & 0.1 & 0.1  & 0.2  & 0.6   & 1.2   & 4.2    & 10.8   \\
\midrule\multirow{3}*{$n=20$}   & \SELRC &0.0 & 0.0 & 0.0  & 0.7  & 4.3   & 54.7  & 1024.0 & \textbf{1800.0} \\
  & \SELQ &0.1 & 0.2 & 0.7  & 1.4  & 5.5   & 12.1  & 43.1   & 149.7  \\
     & \SELQPU & 0.1 & 0.2 & 0.5  & 0.9  & 2.6   & 7.5   & 28.8   & 86.0   \\
\midrule\multirow{3}*{$n=30$}   & \SELRC &0.0 & 0.0 & 0.0  & 1.1  & 8.6   & 169.0 & 1405.1 & \textbf{1800.0} \\
  & \SELQ &0.2 & 0.5 & 2.7  & 7.9  & 15.5  & 43.7  & 159.7  & 475.7  \\
     & \SELQPU & 0.2 & 0.4 & 1.7  & 8.1  & 14.3  & 41.6  & 97.9   & 384.9  \\
\midrule\multirow{3}*{$n=40$}   & \SELRC &0.0 & 0.0 & 0.0  & 1.6  & 11.6  & 228.8 & 1568.8 & \textbf{1800.0} \\
  & \SELQ &1.3 & 1.8 & 5.6  & 28.3 & 35.1  & 135.1 & 361.2  & 1087.7 \\
     & \SELQPU & 1.3 & 1.9 & 3.7  & 14.3 & 46.9  & 147.1 & 310.8  & 958.5  \\
\midrule\multirow{3}*{$n=50$}   & \SELRC &0.0 & 0.0 & 0.0  & 2.1  & 17.0  & 276.6 & 1761.9 & \textbf{1800.0} \\
  & \SELQ &1.4 & 3.0 & 13.9 & 55.8 & 112.6 & 504.9 & 1112.7 & 1580.8 \\
     & \SELQPU & 1.4 & 2.6 & 8.0  & 35.1 & 94.6  & 364.5 & 1056.2 & 1455.1\\\bottomrule

\end{tabular}
\end{table}

\begin{figure}[h!]
\centering
\includegraphics[scale=1]{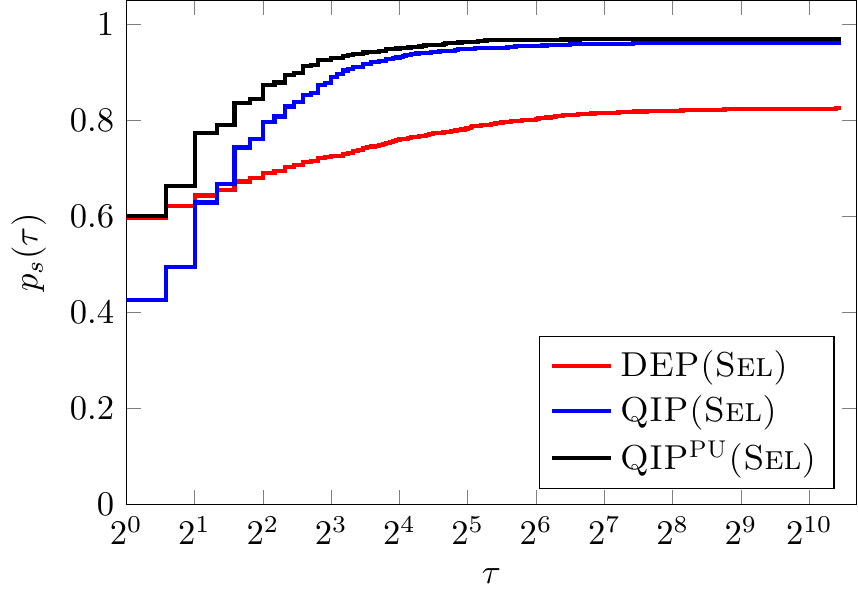}
\caption{Performance profile for all selection instances with $N=4$.}\label{Fig::PerformanceSelectionAll}
\end{figure}
In Figure~\ref{Fig::PerformanceSelectionAll} we provide the performance profile for these experiments.

\clearpage
\section{Additional Results on Multistage Robust Assignment\label{Appendix::Assignment}}

\begin{figure}[h!]
\centering
\includegraphics[scale=1]{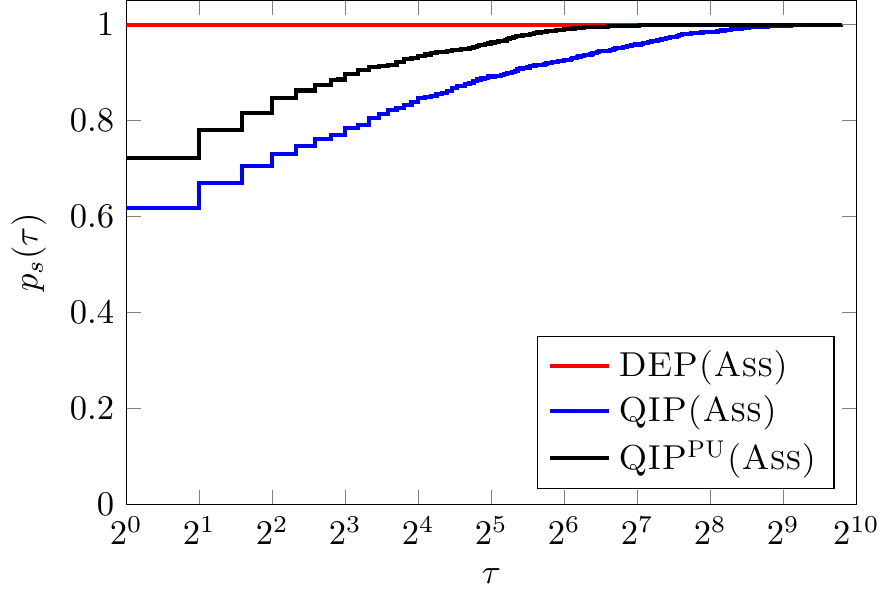}
\caption{Performance profile for all assignment instances with $N=2$.\label{Fig::PerformanceAss_N2}}
\end{figure}

\begin{figure}[h!]
\centering
\includegraphics[scale=1]{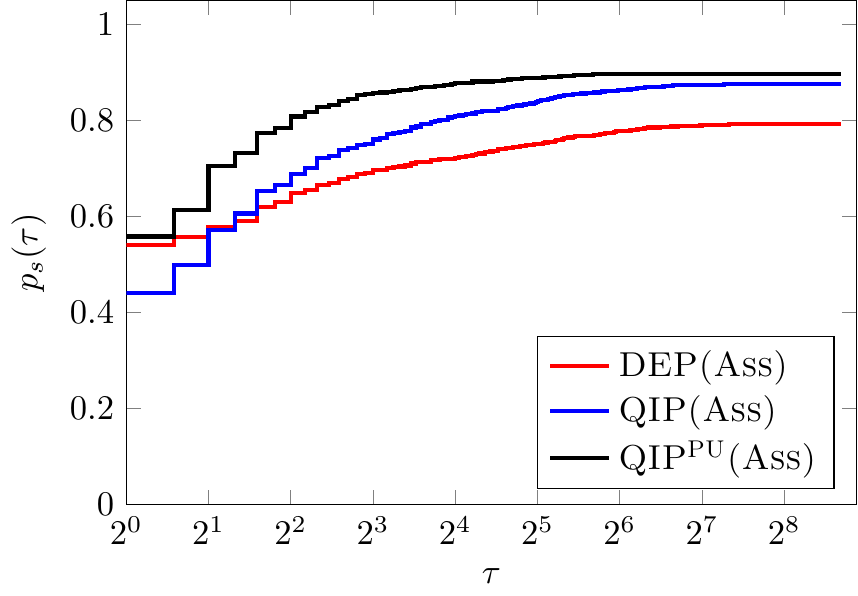}
\caption{Performance profile for all assignment instances with $N=8$.}\label{Fig::PerformanceAss_N8}
\end{figure}

The number of solved instances and the overall average runtime on the assignment problem with $N=2$ and $N=8$ are presented in Table \ref{Table::AssignmentNumTimeN2} and \ref{Table::AssignmentNumTimeN8}, respectively. Figures  \ref{Fig::PerformanceAss_N2} and   \ref{Fig::PerformanceAss_N8} show the corresponding performance profiles.

\begin{table}[h!]
\centering
\caption{Number of solved assignment instances (opt) with $N=2$ and the average runtime (time).}\label{Table::AssignmentNumTimeN2}
\footnotesize
%\begin{tabular}{lp{1.6cm}lp{1.5cm}lp{1.5cm}lp{1.5cm}lp{1.5cm}}
\begin{tabular}{lp{1cm}R{.9cm}rR{.9cm}rR{.9cm}rR{.9cm}r}
\toprule
& & \multicolumn{2}{c}{\hspace{.45cm}$T=1$}    &\multicolumn{2}{c}{\hspace{.45cm}$T=2$}&\multicolumn{2}{c}{\hspace{.45cm}$T=3$}&\multicolumn{2}{c}{\hspace{.45cm}$T=4$} \\
& model& opt &time  & opt &time & opt &time & opt &time     \\
 \midrule \multirow{3}*{$n=4$}  & \ASSRC  & 50     &$ 0.0 $ & 50     &$ 0.0  $ & 50     &$ 0.0  $ & 50     &$ 0.0   $ \\
     & \ASSQ   & 50     &$ 0.0 $ & 50     &$ 0.1  $ & 50     &$ 0.2  $ & 50     &$ 0.2   $ \\
     & \ASSQPU  & 50     &$ 0.1 $ & 50     &$ 0.1  $ & 50     &$ 0.2  $ & 50     &$ 0.2   $ \\
 \midrule \multirow{3}*{$n=5$}  & \ASSRC  & 50     &$ 0.0 $ & 50     &$ 0.0  $ & 50     &$ 0.0  $ & 50     &$ 0.0   $ \\
     & \ASSQ   & 50     &$ 0.1 $ & 50     &$ 0.2  $ & 50     &$ 0.3  $ & 50     &$ 0.7   $ \\
     & \ASSQPU  & 50     &$ 0.1 $ & 50     &$ 0.1  $ & 50     &$ 0.2  $ & 50     &$ 0.4   $ \\
 \midrule \multirow{3}*{$n=6$}  & \ASSRC  & 50     &$ 0.0 $ & 50     &$ 0.0  $ & 50     &$ 0.0  $ & 50     &$ 0.0   $ \\
     & \ASSQ   & 50     &$ 0.1 $ & 50     &$ 0.3  $ & 50     &$ 0.7  $ & 50     &$ 1.7   $ \\
     & \ASSQPU  & 50     &$ 0.1 $ & 50     &$ 0.2  $ & 50     &$ 0.5  $ & 50     &$ 0.8   $ \\
 \midrule \multirow{3}*{$n=7$}  & \ASSRC  & 50     &$ 0.0 $ & 50     &$ 0.0  $ & 50     &$ 0.0  $ & 50     &$ 0.1   $ \\
     & \ASSQ   & 50     &$ 0.2 $ & 50     &$ 0.7  $ & 50     &$ 2.7  $ & 50     &$ 7.6   $ \\
     & \ASSQPU  & 50     &$ 0.1 $ & 50     &$ 0.4  $ & 50     &$ 0.9  $ & 50     &$ 2.1   $ \\
 \midrule \multirow{3}*{$n=8$}  & \ASSRC  & 50     &$ 0.0 $ & 50     &$ 0.0  $ & 50     &$ 0.0  $ & 50     &$ 0.5   $ \\
     & \ASSQ   & 50     &$ 0.2 $ & 50     &$ 2.0  $ & 50     &$ 11.6 $ & 50     &$ 20.1  $ \\
     & \ASSQPU  & 50     &$ 0.2 $ & 50     &$ 0.6  $ & 50     &$ 2.6  $ & 50     &$ 5.1   $ \\
 \midrule \multirow{3}*{$n=9$}  & \ASSRC  & 50     &$ 0.0 $ & 50     &$ 0.0  $ & 50     &$ 0.0  $ & 50     &$ 1.1   $ \\
     & \ASSQ   & 50     &$ 0.2 $ & 50     &$ 5.1  $ & 50     &$ 22.6 $ & 50     &$ 116.3 $ \\
     & \ASSQPU  & 50     &$ 0.2 $ & 50     &$ 1.7  $ & 50     &$ 7.7  $ & 50     &$ 21.4  $ \\
 \midrule \multirow{3}*{$n=10$} & \ASSRC  & 50     &$ 0.0 $ & 50     &$ 0.0  $ & 50     &$ 0.0  $ & 50     &$ 2.1   $ \\
     & \ASSQ   & 50     &$ 0.6 $ & 50     &$ 10.0 $ & 50     &$ 98.0 $ & 50     &$ 351.6 $ \\
     & \ASSQPU  & 50     &$ 0.3 $ & 50     &$ 4.6  $ & 50     &$ 22.9 $ & 50     &$ 84.3  $ \\\bottomrule
\end{tabular}
\end{table}

\begin{table}[h!]
\centering
\vspace{-.6cm}~
\caption{Number of solved assignment instances (opt) with $N=8$ and the average runtime (time).}\label{Table::AssignmentNumTimeN8}
\footnotesize
\begin{tabular}{lp{1cm}R{.9cm}rR{.9cm}rR{.9cm}rR{.9cm}r}

\toprule
& & \multicolumn{2}{c}{\hspace{.45cm}$T=1$}    &\multicolumn{2}{c}{\hspace{.45cm}$T=2$}&\multicolumn{2}{c}{\hspace{.45cm}$T=3$}&\multicolumn{2}{c}{\hspace{.45cm}$T=4$}    \\
& model& opt &time  & opt &time & opt &time & opt &time     \\
 \midrule \multirow{3}*{$n=4$}  & \ASSRC & 50 &$ 0.0 $ & 50 &$ 0.0   $      & 50 &$ 3.2    $ & 50 &$ 118.6  $ \\
     & \ASSQ  & 50 &$ 0.1 $ & 50 &$ 0.2   $                                 & 50 &$ 0.8    $ & 50 &$ 4.0    $ \\
     & \ASSQPU & 50 &$ 0.1 $ & 50 &$ 0.2   $                                & 50 &$ 1.3    $ & 50 &$ 5.8    $ \\
 \midrule \multirow{3}*{$n=5$}    & \ASSRC & 50 &$ 0.0 $ & 50 &$ 0.2   $    & 50 &$ 11.8   $ & 43 &$ 552.1  $ \\
     & \ASSQ  & 50 &$ 0.1 $ & 50 &$ 0.6   $                                 & 50 &$ 2.5    $ & 50 &$ 13.9   $ \\
     & \ASSQPU & 50 &$ 0.2 $ & 50 &$ 0.4   $                                & 50 &$ 2.4    $ & 50 &$ 18.4   $ \\
 \midrule \multirow{3}*{$n=6$}   & \ASSRC & 50 &$ 0.0 $ & 50 &$ 0.9   $     & 50 &$ 70.0   $ & 25 &$ 1231.4  $ \\
     & \ASSQ  & 50 &$ 0.3 $ & 50 &$ 2.5   $                                 & 50 &$ 13.6   $ & 50 &$ 63.6   $ \\
     & \ASSQPU & 50 &$ 0.2 $ & 50 &$ 1.7   $                                & 50 &$ 8.4    $ & 50 &$ 70.0   $ \\
 \midrule \multirow{3}*{$n=7$}   & \ASSRC & 50 &$ 0.0 $ & 50 &$ 1.4   $     & 49 &$ 246.9  $ & 7  &$ 1654.8 $ \\
     & \ASSQ  & 50 &$ 0.4 $ & 50 &$ 11.3  $                                 & 50 &$ 58.8   $ & 49 &$ 340.5  $ \\
     & \ASSQPU & 50 &$ 0.3 $ & 50 &$ 4.1   $                                & 50 &$ 33.7   $ & 50 &$ 248.0  $ \\
 \midrule \multirow{3}*{$n=8$}   & \ASSRC & 50 &$ 0.0 $ & 50 &$ 4.8   $     & 41 &$ 689.6  $ & 2  &$ 1747.9  $ \\
     & \ASSQ  & 50 &$ 1.2 $ & 50 &$ 47.4  $                                 & 50 &$ 251.1  $ & 39 &$ 1062.5  $ \\
     & \ASSQPU & 50 &$ 0.5 $ & 50 &$ 13.5  $                                & 50 &$ 110.8  $ & 44 &$ 981.1  $ \\
 \midrule \multirow{3}*{$n=9$}    & \ASSRC & 50 &$ 0.0 $ & 50 &$ 21.9  $    & 24 &$ 1161.9  $ & 3  &$ 1751.5  $ \\
     & \ASSQ  & 50 &$ 2.9 $ & 50 &$ 228.8 $                                 & 31 &$ 1377.1 $ & 7  &$ 1673.8  $ \\
     & \ASSQPU & 50 &$ 1.5 $ & 50 &$ 45.8  $                                & 50 &$ 598.5  $ & 4  &$ 1746.5 $ \\
 \midrule \multirow{3}*{$n=10$}   & \ASSRC & 50 &$ 0.0 $ & 50 &$ 56.4  $    & 15 &$ 1415.2  $ & 1  &$ 1771.8  $ \\
     & \ASSQ  & 50 &$ 8.1 $ & 43 &$ 797.2 $                                 & 7  &$ 1642.6  $ & 1  &$ 1767.5  $ \\
     & \ASSQPU & 50 &$ 3.0 $ & 50 &$ 266.5 $                                & 8  &$ 1680.9 $ & 1  &$ 1787.1 $\\\bottomrule
\end{tabular}
\end{table}
\end{document}